\documentclass[12pt]{amsart}
\usepackage{latexsym,amssymb,amsmath,amsthm}
\usepackage{comment}
\usepackage{float}
\usepackage{graphicx}
\usepackage{color}
\usepackage{listings}
\definecolor{color7}{rgb}{0.8,0.8,1}
\definecolor{color8}{rgb}{1,1,0.6}
\theoremstyle{definition}

\newcounter{slide} 
\def\slide#1#2{ \item 
\vspace{5mm}
\parbox{14.5cm}{
\parbox{6.4cm}{\scalebox{0.30}{\includegraphics{slides/#1.pdf}}}
\hspace{5mm}
\parbox{7.0cm}{{\bf #2}}
} 
\vspace{5mm}
} 

\newcounter{problem}
\def\problem#1{ \item #1 \vspace{3mm} }
\newcounter{problemm}
\def\problemm#1{ \item #1 \vspace{3mm} }

\def\G{\mathcal{G}}

\def\Sin{\rm{Sin}}

\title{If Archimedes would have known functions ...}
\author{Oliver Knill}
\date{March 23, 2014}
\address{
        Department of Mathematics \\
        Harvard University \\
        Cambridge, MA, 02138
        }

\begin{document}
\maketitle

\begin{abstract}
Could calculus on graphs have emerged by the time of Archimedes,
if function, graph theory and matrix concepts were available 2300 years ago?
\end{abstract}

\section{Single variable calculus} 

Calculus on integers deals with {\bf functions}
$f(x)$ like $f(x)=x^2$. The {\bf difference} $Df(x)=f'(x) = f(x+1)-f(x) = 2x+1$ as well as
the {\bf sum} $Sf(x)=\sum_{k=0}^{x-1} f(k)$ with the understanding $Sf(0)=0,Sf(-1)=f(-1)$ are
functions again. We call $Df$ the {\bf derivative} and $Sf$ the {\bf integral}. The identities
$DSf(x)=f(x)$ and $SDf(x)=f(x)-f(0)$ are the {\bf fundamental theorem of calculus}.
Linking sums and differences allows to compute sums (which is difficult
in general) by studying differences (which is easy in general). Studying derivatives of
basic functions like $x^n,\exp(a \cdot x)$ will allow to sum such functions. 
As operators, $Xf(x) = x s^* f(x)$ and $Df(x)=[s,f]$ where $s f(x)=f(x+1)$ and
$s^* f(x) = f(x-1)$ are translations. We have $1 x = x, X x = x (x-1), X^2 x = x (x-1)(x-2)$. 
The derivative  operator $Df(x) = (f(x+1)-f(x))s$ satisfies the {\bf Leibniz product rule} 
$D(f g) = (Df) g + (Dg) f = Df g + f^+ Dg$.
Since {\bf momentum} $P=iD$ satisfies the {\bf anti-commutation relation} $[X,P]=i$ we have
{\bf quantum calculus}. The polynomials $[x]^n = X^n x$ satisfy $D[x]^n=n [x]^{n-1}$ and 
$X^k X^m = X^{k+m}$. The {\bf exponential} $\exp(a \cdot x) = (1+a)^x$
satisfies $D\exp(a \cdot x) = a \exp(a \cdot x)$ and $\exp(a \cdot (x+y)) = \exp(a \cdot x)  \exp(a \cdot y)$.  
Define $\sin(a \cdot x)$ and $\cos(a \cdot x)$ as the real and imaginary part of $\exp(i a \cdot x)$. 
From $D e^{a \cdot x}=a e^{a \cdot x}$, we deduce $D\sin(a \cdot x) = a \cos(a \cdot x)$ and 
$D \cos(a \cdot x) = -a \sin(a \cdot x)$. 
Since $\exp$ is monotone, its inverse $\log$ is defined. Define the reciprocal $X^{-1}$ by 
$X^{-1} x = D \log(x)$ and $X^{-k} = (X^{-1})^k$. We check identities like 
$\log(x y) = \log(x) + \log(y)$ from which $\log(1)=0$ follow. 
The {\bf Taylor formula} $f(x) = \sum_{k=0}^{\infty} D^kf(0) X^k/k!$ is
due to {\bf Newton and Gregory} and can be shown by induction or by noting that $f(x+t)$ satisfies
the transport equation $D_t f = D f$ which is solved by $f(t) = \exp(D t) f(0)$.
The Taylor expansion gives $\exp(x) = \sum_{k=0}^{\infty} X^k/k!$.
We also get from $(1-x)^{-1} = \sum_{k=0}^{\infty} X^k$ that $\log(1-x) = -\sum_{k=1}^{\infty} X^k/k$.
The functions $c \sin(a \cdot x) + d \cos(a \cdot x)$ solve the {\bf harmonic oscillator} equation
$D^2 f = -a^2 f$. The equation $Df = g$ is solved by the {\bf anti derivative} $f = S g(x) + c$.
The {\bf wave equation} $D_t^2 f = D^2 f$ can be written as $(D_t - D) (D_t + D) f=0$
which thus has the solution $f(t)$ given as a linear combination of 
$\exp(D t) f(x) = f(x+t)$ and $\exp(-Dt) f(x) = f(x-t)$.
We replace now integers with a graph, the derivative $iD$ with a matrix $D$,
functions by {\bf forms}, anti-derivatives by $D^{-1} f$, and Taylor with Feynman.

\section{Multivariable calculus} 

{\bf Space} is a finite simple graph $G=(V,E)$.
A subgraph $H$ of $G$ is a {\bf curve} if every {\bf unit sphere} of a vertex $v$ in $H$ 
consists of $1$ or $2$ disconnected points. A subgraph $H$ is a {\bf surface } if every unit sphere is 
either an interval graph or a circular graph $C_n$ with $n \geq 4$. Let $\G_k$ denote the set of $K_{k+1}$ 
subgraphs of $G$ so that $G_0=V$ and $G_1=E$ and $\G_2$ is the set of triangles in $G$. Intersection relations
render each $\G_k$ a graph too. The linear space $\Omega_k$ 
of anti-symmetric functions on $\G_k$ is the set of {\bf $k$-forms}. For $k=0$, we have 
{\bf scalar functions}. To fix a basis, equip simplices like edges and triangles with orientations.
Let $d$ denote the exterior derivative $\Omega_k \to \Omega_{k+1}$. If $\delta H$ be the {\bf boundary} of $H$,
then $df(x)=f(\delta x)$ for any $x \in \G_k$. The boundary of a curve $C$ is empty or consists of two points. 
The boundary of a surface $S$ is a finite set of closed curves. 
A closed surface $S$ is a discretisation of a $2$-dimensional classical surface as we can build an atlas from wheel subgraphs.
An {\bf orientation} of $S$ a choice of orientations on triangles in $S$ such that the boundary 
orientation of two intersecting triangles match. If $F$ is a $2$-form and $S$ is an orientable
surface, define $\int_S F \; dS =\sum_{x \in \G_2(S)} F(x)$.
An orientation of a curve $C$ is a choice of orientations on edges in $C$ such that the boundary 
orientation of two adjacent edges correspond. If $f$ is a $1$-form and $C$ is orientable, 
denote by $\int_C F \; ds$  the sum $\sum_{x \in G_1(C)} F(x)$. For a $0$-form, write
${\rm grad}(f) = df$ and for a $1$-form $F$, write ${\rm curl}(F)=dF$. 
If $S$ is a surface and $F$ is a $1$-form, the {\bf theorem of Stokes} tells that 
$\int_S dF = \int_{\delta S} F$. It is proved by induction using that
for a single simplex, it is the definition. If $H$ is a $1$-dimensional graph, and $F$ is a function, 
then $\int_H dF = \int_{\delta H} F$ is the {\bf fundamental theorem of line integrals}. 
A graph is a {\bf solid}, if every unit sphere is either a triangularization of the
$2$-sphere or a triangularization of a disc. Points for which the sphere is a disc form the
{\bf boundary}. If $F$ is a $2$-form and $H$ is a solid, then $\int_H dF = \int_{\delta H} F $
is {\bf Gauss theorem}. The divergence is the adjoint $d^*: \Omega_1 \to \Omega_0$. In school calculus, 
$1$-forms and $2$-forms are often identified as {\bf ``vector fields"} and $3$-forms on the solid $H$ 
is treated as {\bf scalar functions}. Gauss is now the {\bf divergence theorem}. Define the {\bf Dirac matrix} $D=d+d^*$
and the {\bf form Laplacian} $L=d d^* + d^* d$. Restricted to scalar functions, it is $d^* d f= {\rm div}({\rm grad}f)$
which is equal to $A-B$, where $A$ is the adjacency matrix and $B$ the diagonal degree matrix of the graph. 
The {\bf Schr\"odinger equation} $f' = iD f$ is solved by
$f(t) = e^{it D }f(0) = f(0) + t D f(0) + D^2 t^2 f(0)/2! + \dots$. It is a
{\bf multi-dimensional Taylor equation}, because by matrix multiplication, the term $(t^n D^n) f$ 
is the sum over all possible paths $\gamma$ of length $n$ where each step $k$ is multiplied with $t f(\gamma(k))$.
The solution of the Schr\"odinger equation $f(t,y)$ is now the average over all possible paths from $x$ to $y$.
This {\bf Feynman-Kac} formula holds for general $k$-forms. Given an initial form $f$, it describes $f(t)$. 
The {\bf heat equation} $f'=-L f$  has the solution $e^{-Lt}f(0)$, the heat flow. Because $L$ leaves $k$-forms invariant,
the Feynman-Kac formula deals with paths which stay on the graph $\G_k$. 
As any linear differential equation, it can be solved using eigenfunctions of $L$. If $G=C_n$, it is a discrete Fourier basis. 
The wave equation $f''=-Lf$ can be written as $(\partial_t - i D) (\partial_t + i D)f =0$
which shows that $f =\cos(D t) f(0) + \sin(D t) D^{-1} f'(0) = e^{i D t} \psi$ with $\psi = f(0) - i D^{-1} f'(0)$. 
Feynman-Kac for a $k$-form expresses $\psi(t,x)$ as the expectation ${\rm E}[\sum_{x \in \gamma} t^{n(\gamma)} f(0,x)]$ on 
a probability space of finite paths $\gamma$ in the simplex graph $\G_k$
where $n(\gamma)$ is the length of $\gamma$. 

\pagebreak

\section{Pecha-Kucha}

{\bf Pecha-Kucha} is a presentation form in which 20 slides are shown exactly 20 seconds each. 
Here are slides for a talk given on March 6, 2013. The bold text had been prepared to be spoken,
the text after consists of remarks posted on March 2013.
\begin{list}{\fcolorbox{color7}{color7}{\Large{\arabic{slide}}}}{\usecounter{slide}}

\slide{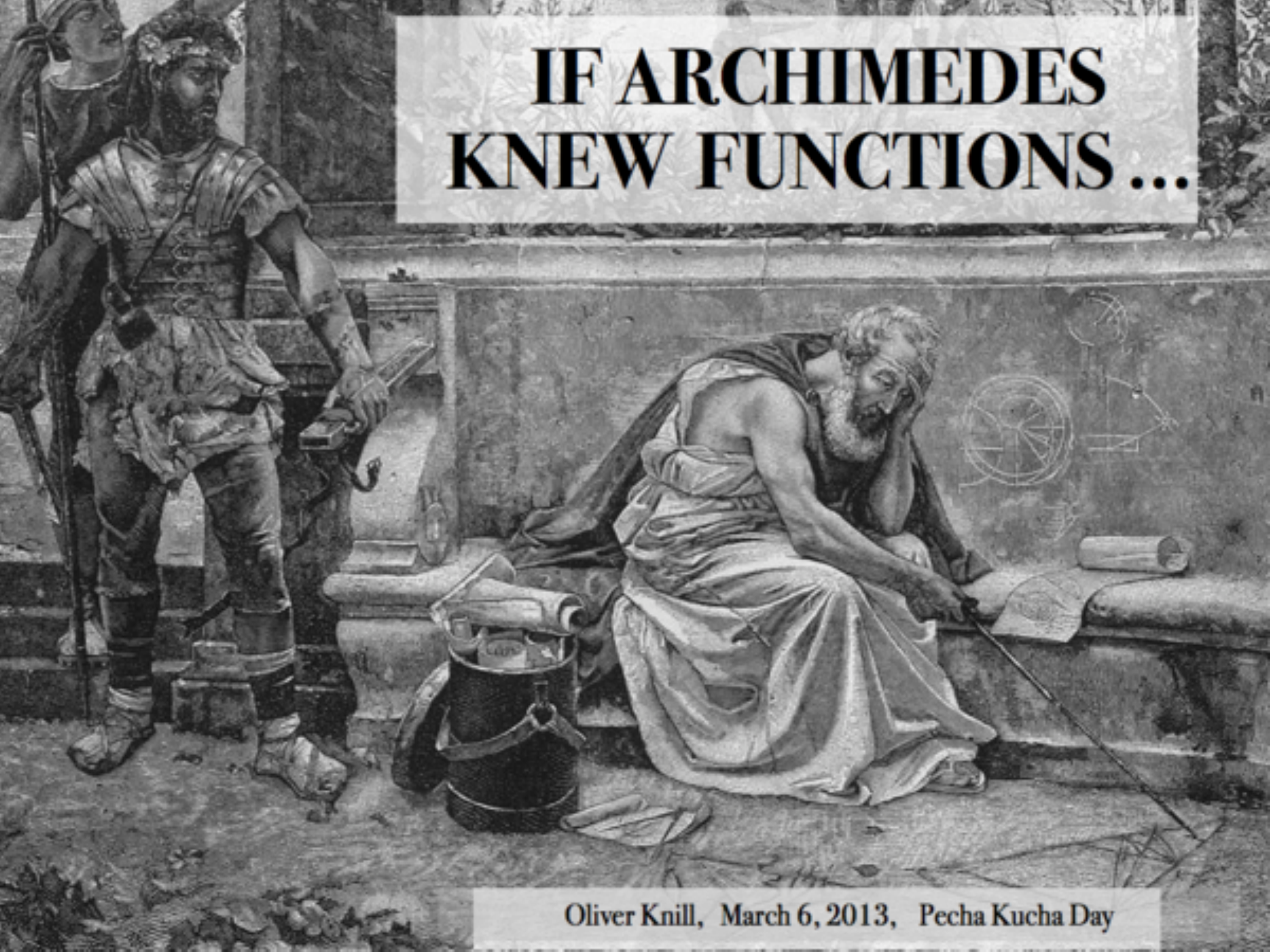}{What if Archimedes would have known the concept of a function? The following story
borrows from classes taught to students in singular variable calculus and at the extension school.
The mathematics is not new but illustrates part of what one calls ``quantum calculus" or ``calculus
without limits".}

Quantum calculus comes in different flavors. A book with that title was written by Kac and Cheung
\cite{KacCheung}. 
The field has connections with numerical analysis, combinatorics and even number theory. The calculus of finite
differences was formalized by George Boole but is much older.
Even Archimedes and Euler thought about calculus that way. It also has connections
with nonstandard analysis. In "Math 1a, a single variable calculus course 
or Math E320, a Harvard extension school course, I spend 2 hours each semester with such material. 
Also, for the last 10 years, a small slide show at the end of the semester in multivariable calculus 
featured a presentation "Beyond calculus". 

\slide{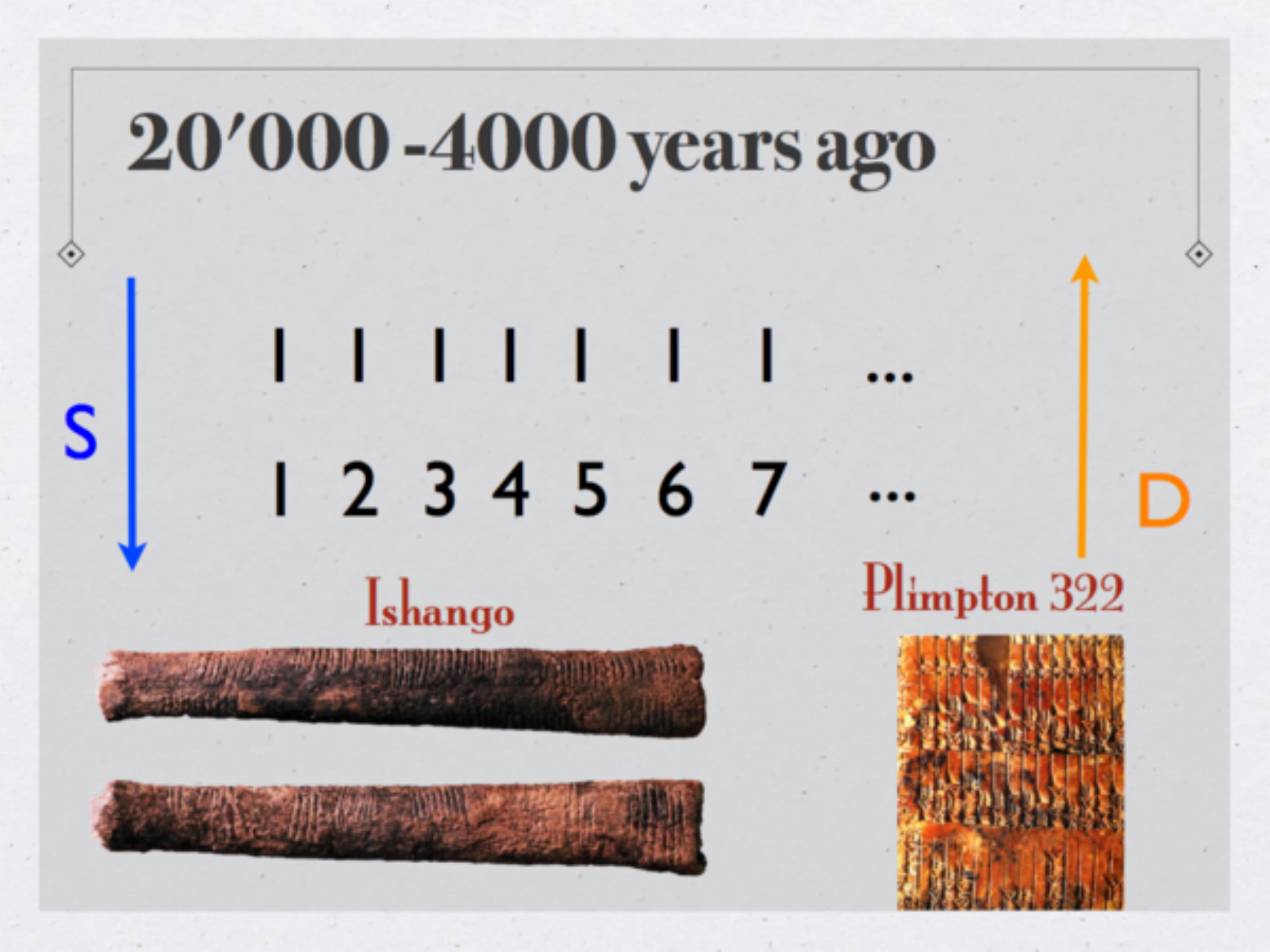}{About 20 thousand years ago, mathematicians represented numbers as marks on bones.
I want you to see this as the constant function 1.
Summing up the function gives the concept of number. 
The number 4 for example is $1+1+1+1$. Taking differences brings us back: $4-3=1$.}

The {\bf Ishango bone} displayed in the slide
is now believed to be about 20'000 years old. The number systems
have appeared in the next couple of thousand years, independently in different places. 
They were definitely in place 4000 BC because we have Carbon dated Clay tablets featuring numbers. 
Its fun to make your own Clay tablets: either on clay or
on chewing gum.

\slide{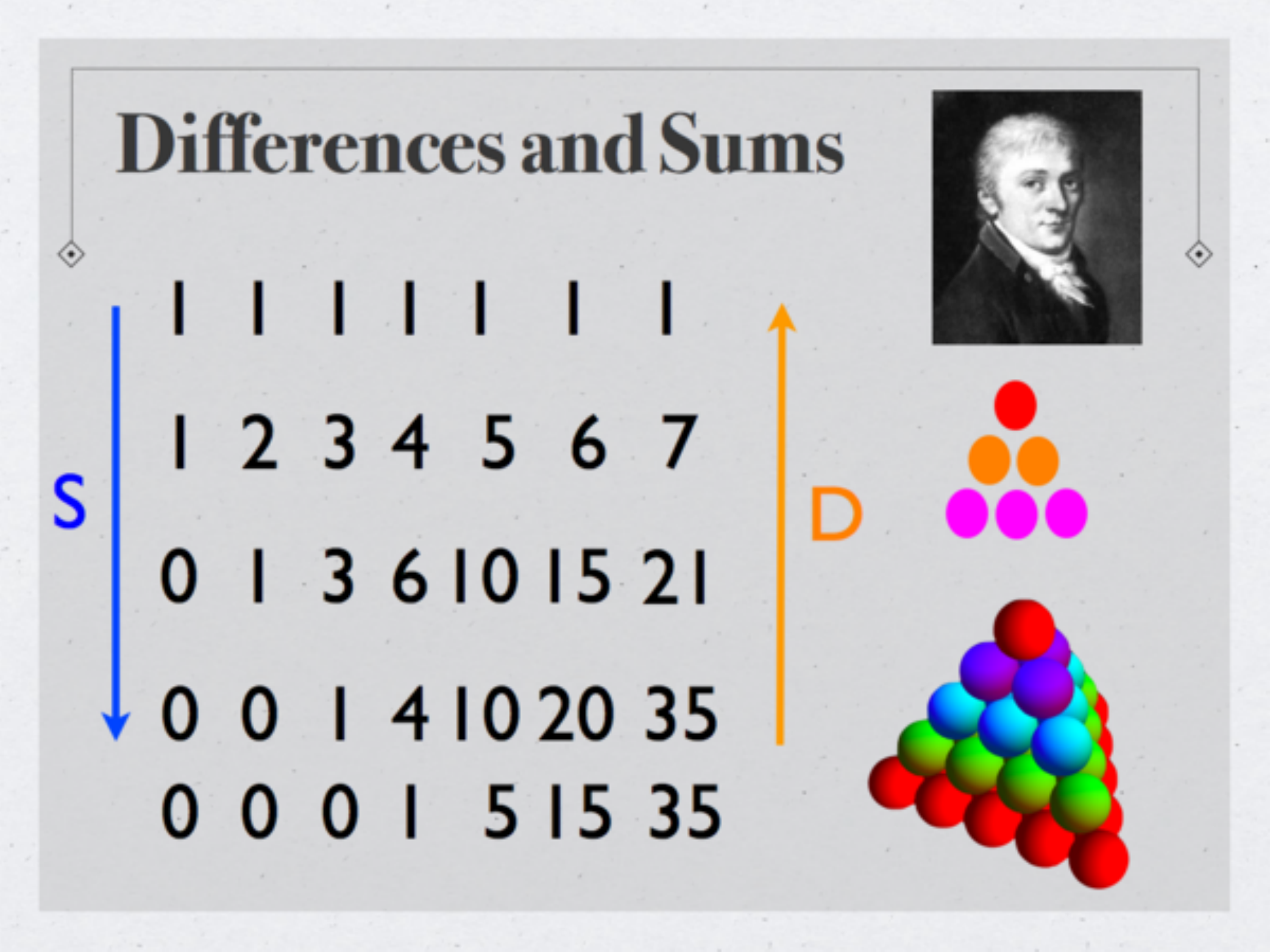}{
When summing up the previous function $f(x)=x$, we get {\bf triangular numbers}.
They represent the area of triangles.  Gauss got a formula for them as a school kid. 
Summing them up gives tetrahedral numbers, which represent volumes. Integration is 
summation, differentiation is taking differences.  }

The numbers $[x]^n/n!$ are examples of {\bf polytopic numbers}, numbers which 
represent patterns. Examples are triangular, tetrahedral or pentatopic numbers. Since we use the forward
difference $f(x+1)-f(x)$ and start summing from $0$, our formulas are shifted. For example, the triangular numbers
are traditionally labeled $n(n+1)/2$, while we write $n(n-1)/2$. 
We can think about the new functions $[x]^n$ as a 
new basis in the linear space of polynomials. 

\slide{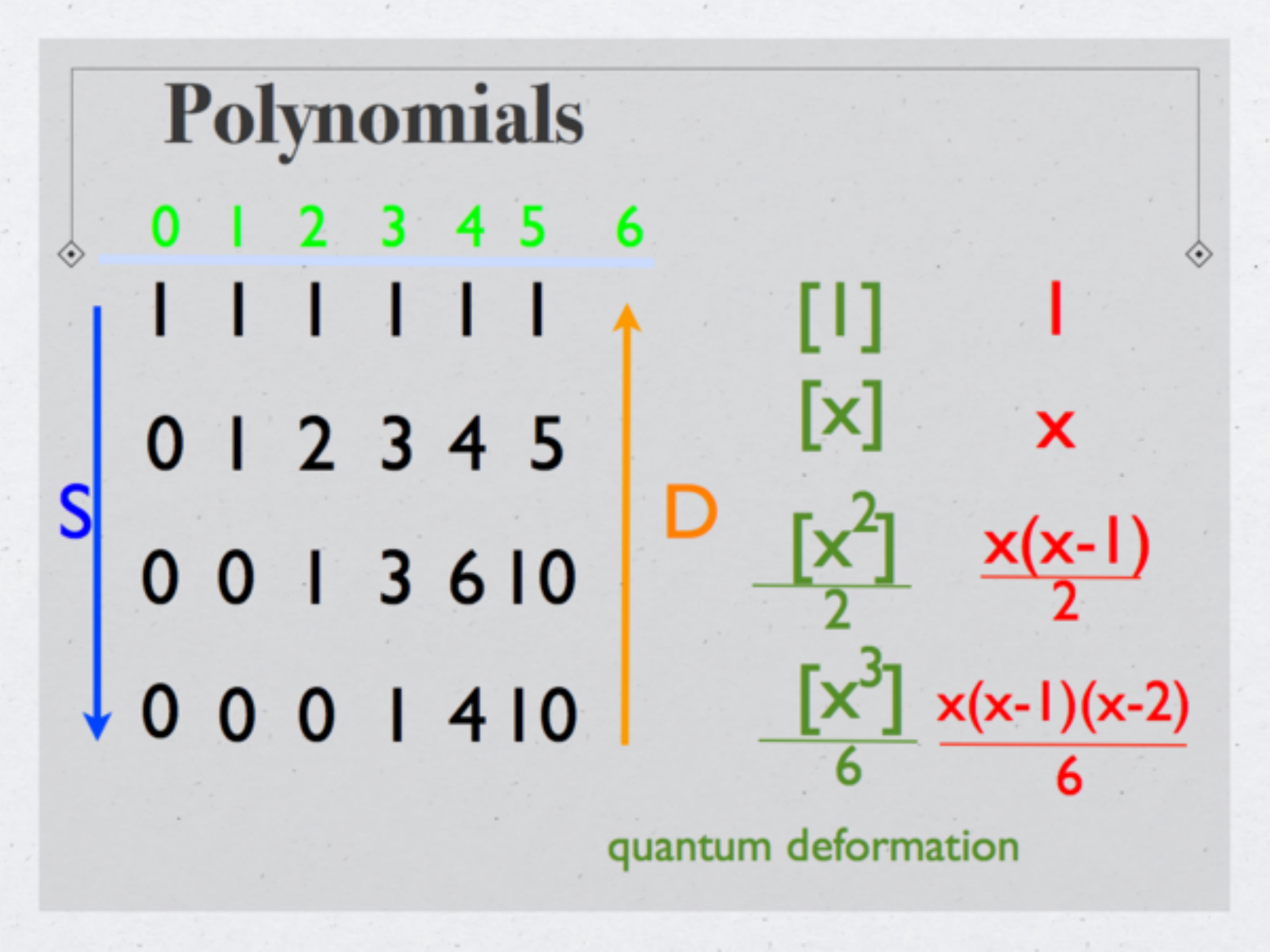}{
The polynomials appear as part of the Pascal triangle.
We see that the summation process is the inverse of the difference process. 
We use the notation $[x]^n/n!$ for the function in the nth row.}

The ``renaming idea" $[x]^n$ is part of {\bf quantum calculus}. 
It is more natural when formulated in a non-commutative algebra, a crossed product of
the commutative algebra we are familiar with. The commutative $X=C^*$-algebra $C(R)$ of all continuous real functions encodes
the topology of the real numbers. If $s(x)=x+1$ is translation, we can look at the algebra of operators on 
$H=L^2(R)$ generated by $X$ and the translation operator s. With  $[x]= x s^*$, the operators $[x]^n$ 
are the multiplication operators generated by the polynomials of the slide. The derivative $Df$
be written as the commutator $Df = [s,f]$. In quantum mechanics, functions are
treated as operators. The deformed algebra is no more commutative: with $Qf(x)=x s^*$ satisfying $Q^n=[x^n]$
and $Pf=iDf$ the anti-commutation relations $[Q,P]=QP-PQ=i$ hold.
Hence the name ``quantum calculus".

\slide{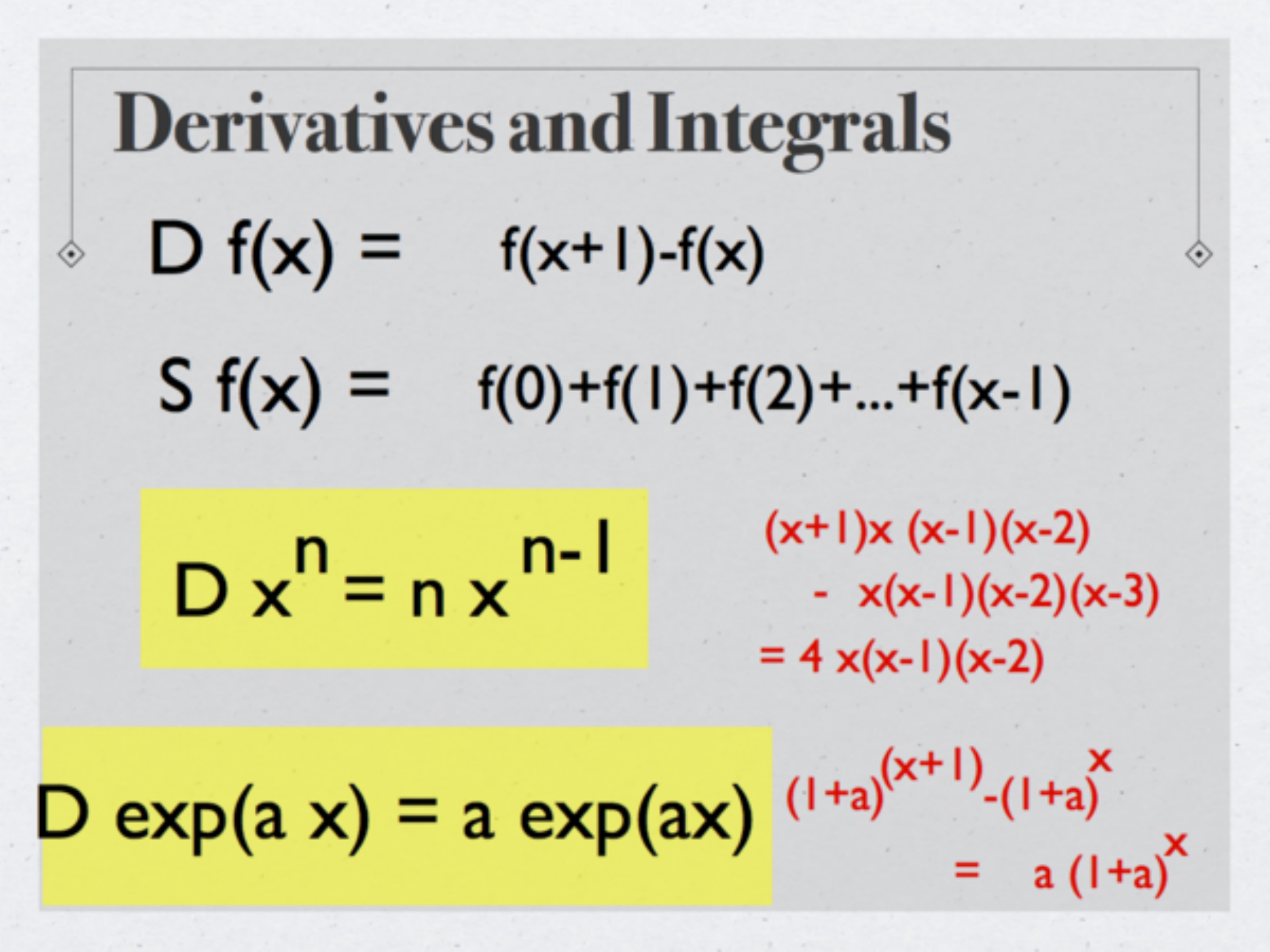}{
With the difference operation $D$ and the summation operation $S$, we already get an important formula. 
It tells that the derivative of $x^n$  is $n$ times x to the $n-1$.
Lets also introduce the function $\exp(a x) = (1+a)$ to the power $x$. This is the 
compound interest formula. We check that its derivative is a constant times the function itself.
}

The deformed exponential is just a rescaled exponential with base $(1+a)$.  Writing it as a compound
interest formula allows to see the formula $\exp'(ax)=a \exp(ax)$ as the property
that the bank pays you a multiple a $\exp(ax)$ to your fortune $\exp(ax)$. The compound interest formula appears in 
the movie "The bank".
We will below look at the more general exponential $\exp_h(x) = (1+a h)^{x/h}$ which is the exponential 
$\exp(a x)$ to the "Planck constant" $h$. If $h$ goes to zero, we are lead to the standard 
exponential function $\exp(x)$. One can establish the limit because 
$\exp_h(a x) \le \exp(a x)\le (1+a h) \exp_h(a x)$.  

\slide{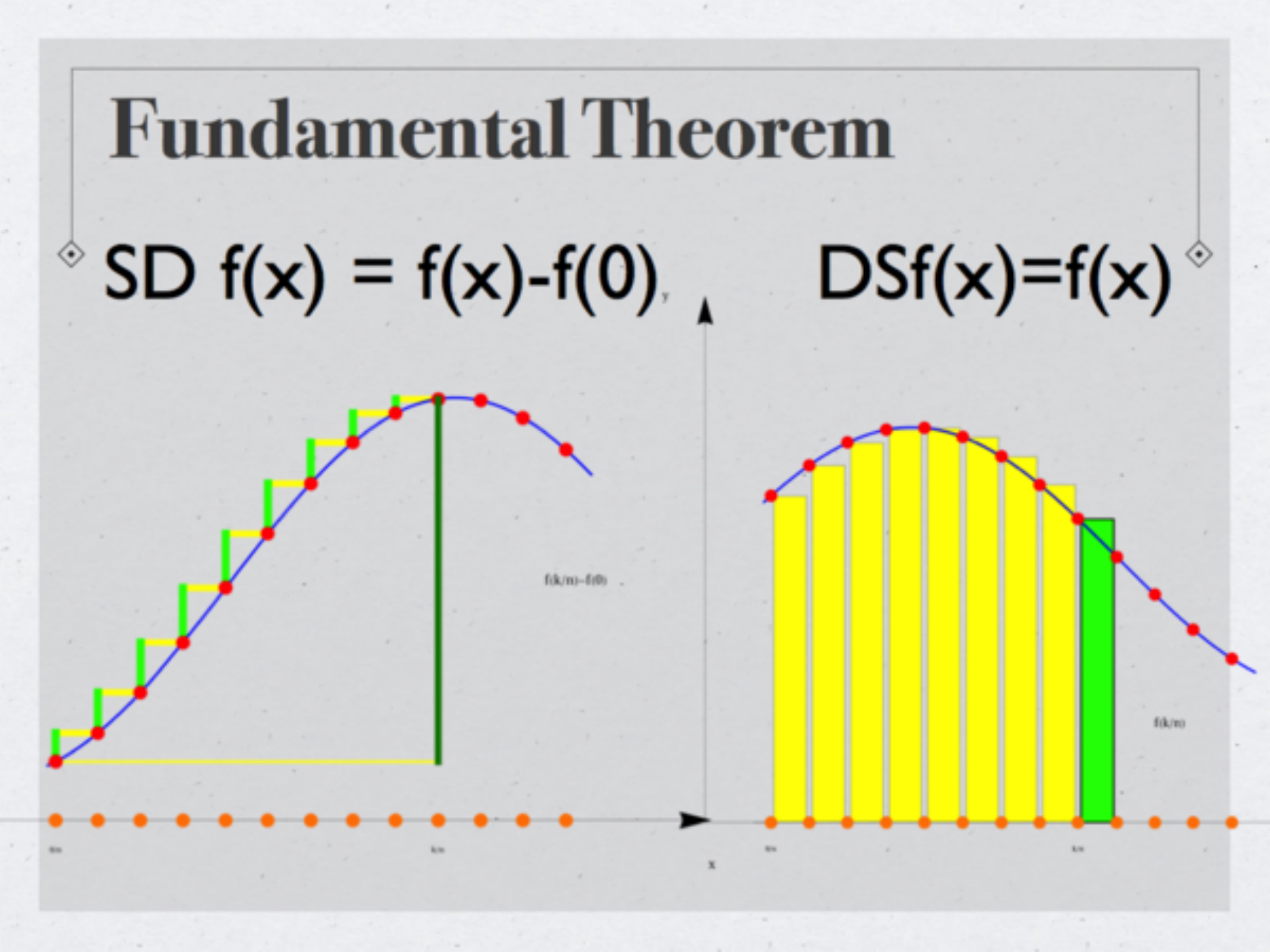}{
The heart of calculus is the fundamental theorem of calculus. If we take the sum of the
differences, we get the difference of f(x) minus f(0). If we take the difference of the sum, we get
the function f(x). The pictures pretty much prove this without words. These formulas are true for any function.
They do not even have to be continuous.}

Even so the notation $D$ is simpler than $d/dx$ and $S$ is simpler than the integral sign, it is the language with 
$D,S$ which makes the theorem look unfamiliar. Notation is very important in mathematics. 
Unfamiliar notation can lead to rejection. One of the main points to be made
here is that we do not have to change the language. We can leave all calculus books as they are and
just look at the content in new eyes. The talk is advertisement for the calculus we know, teach
and cherish. It is a much richer theory than anticipated.

\slide{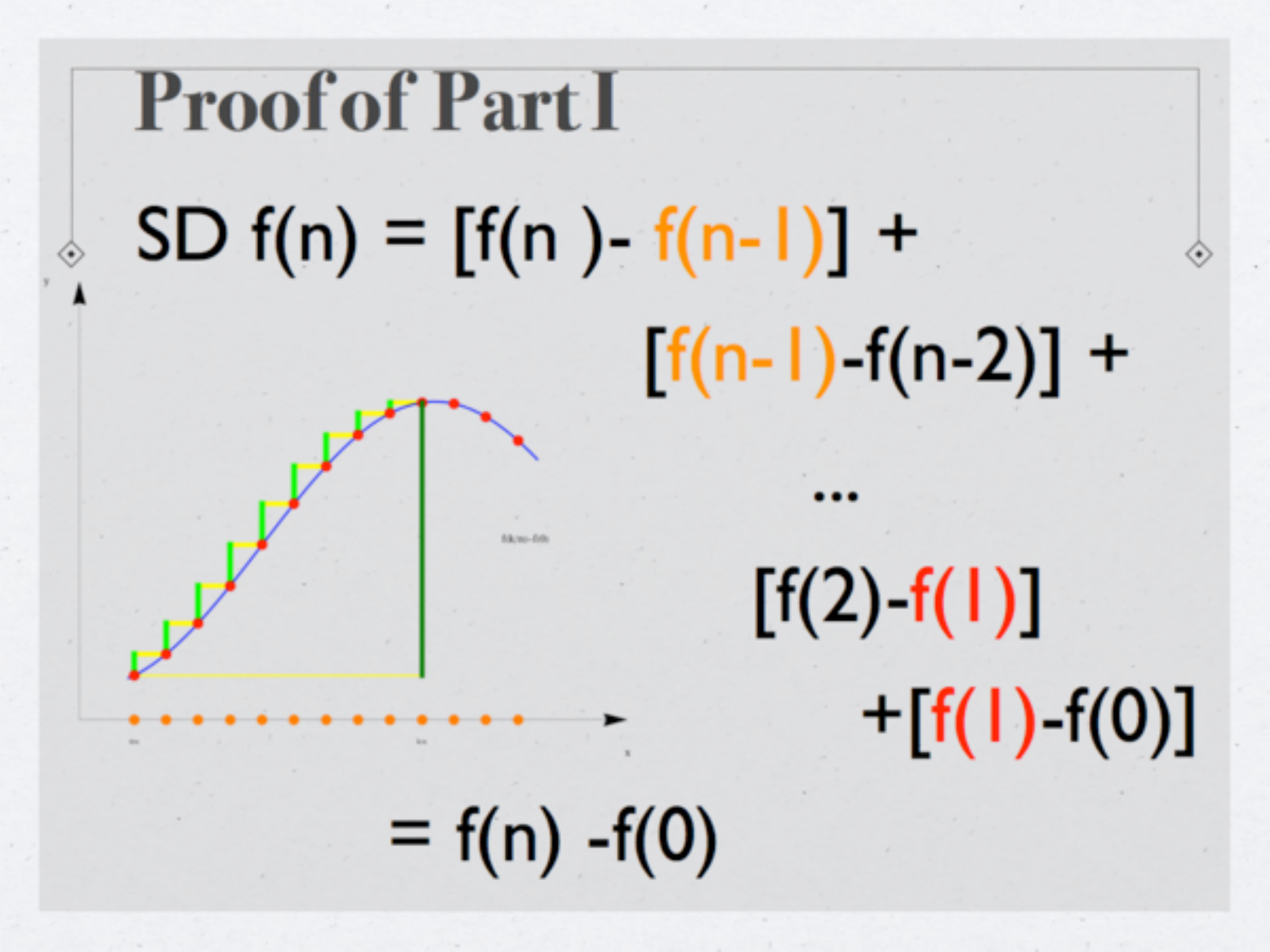}{
The proof of part I shows a telescopic sum. The cancellations
are at the core of all fundamental theorems. They appear also in multivariable calculus 
in differential geometry and beyond.}

It is refreshing that one can show the proof of the fundamental theorem of calculus
early in a course. Traditionally, it takes much longer until one reaches the point.
Students have here only been puzzled about the fact that the result holds only for x=nh and not in general.
It would already here mathematically make more sense to talk about a result on a finite linear graph but that
would increase mentally the distance to the traditional calculus in
\cite{Bressoud} on "Historical Reflections on Teaching the fundamental Theorem of Integral Calculus". 
Bressoud tells that "there is a fundamental problem with the statement of the FTC and that only a 
few students understand it".  
I hope this Pecha-Kucha can help to see this important theorem from a different perspective. 

\slide{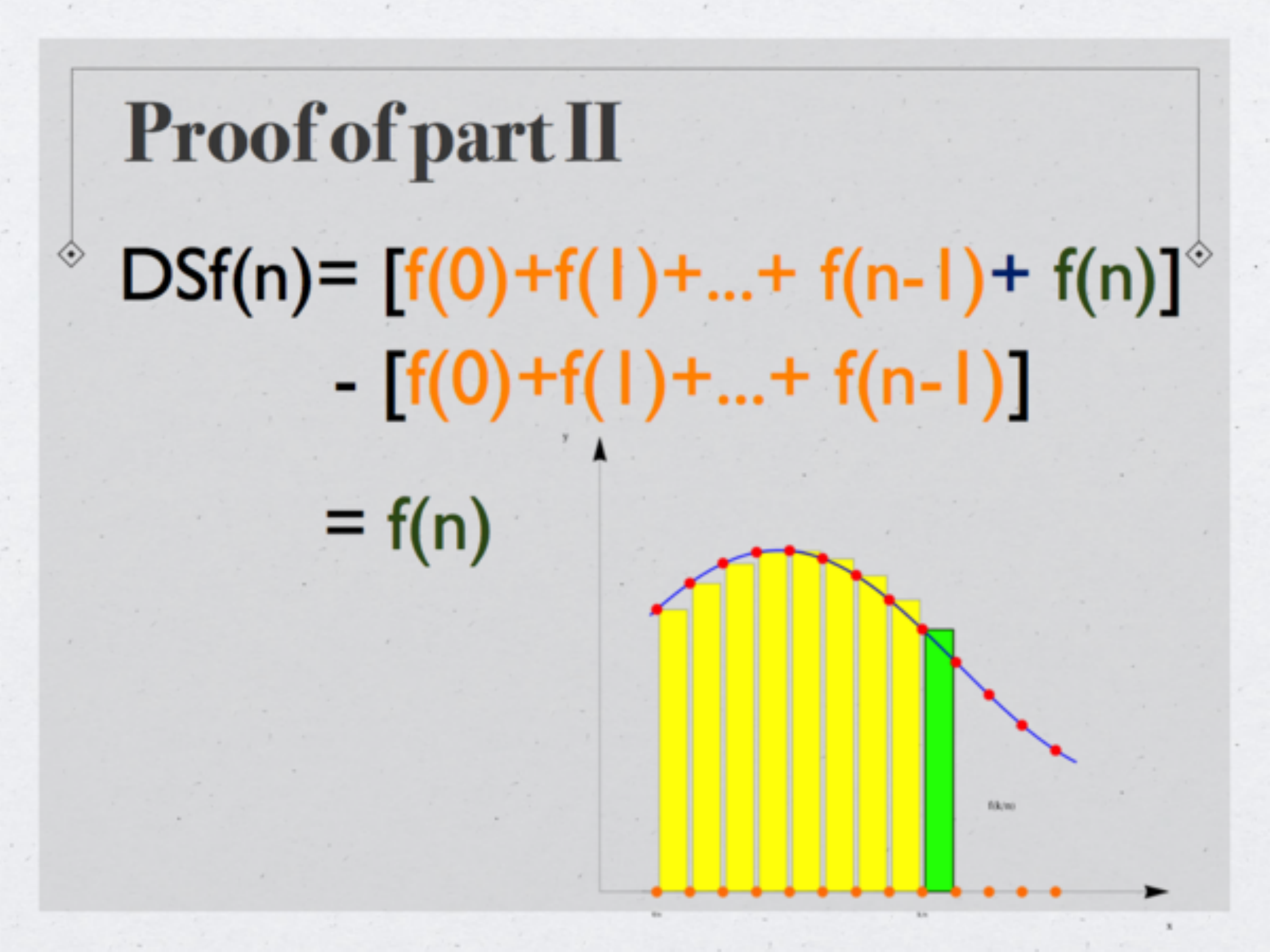}{
The proof of the second part is even simpler. The entire proof can be redone, when the step size h=1 is
replaced by a positive h. The fundamental theorem allows us to solve a difficult problem:
summation is global and hard, differences is local and  easy. By combining the two we can make the hard
stuff easy. }

This is the main point of calculus. Integration is difficult, Differentiation is easy. Having them 
linked makes the hard things easy. At least part of it. 
It is the main reason why calculus is so effective. These ideas
go over to the discrete. Seeing the fundamental idea first in a discrete setup can help
to take the limit when h goes to zero. 

\slide{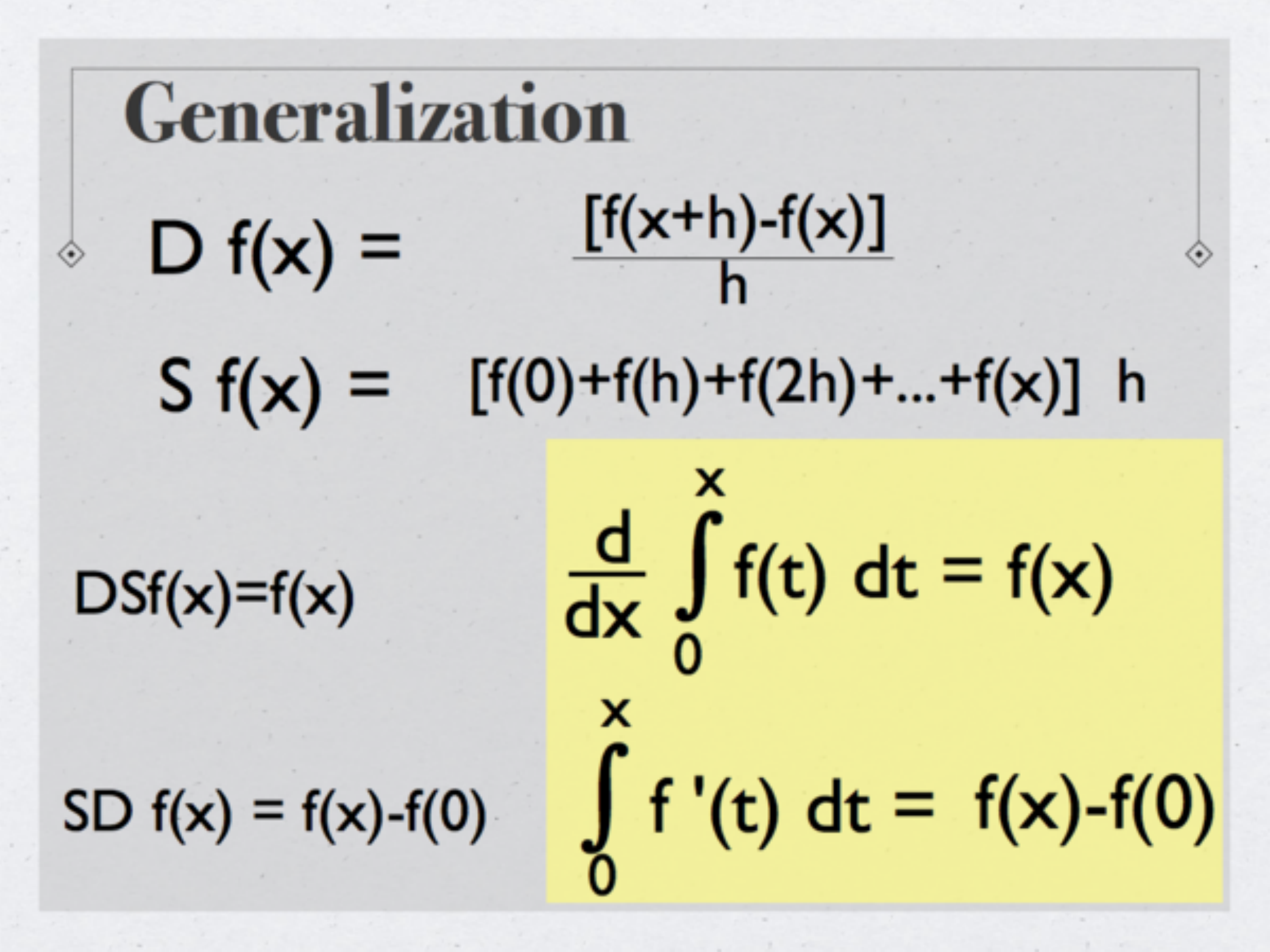}{
We can adapt the step size h from 1 to become any positive number h.
The traditional fundamental theorem of calculus is stated here also in the notation of Leibniz.
In classical calculus we would take the limit h to 0, but we do not go that way in this talk.}

When teaching this, students feel a slight unease with the additional variable h. 
For mathematicians it is natural to have free parameters, for students who have had little
exposure to mathematics, it is difficult to distinguish between the variable x and h. 
This is the main reason why in this talk, I mostly stuck to the case h=1. 
In some sense, the limit h to 0 is an idealization. This is why the limit h going to 0
is so desirable. But we pay a prize: the class of function we can deal with, is much smaller. 
I myself like to see traditional calculus as a limiting case of a larger theory. The
limit h to zero is elegant and clean. But it takes a considerable effort at first 
to learn what a limit is. Everybody who teaches the subject can confirm the principle
that often things which were historically hard to find, are also harder to teach and master.
Nature might have given up on it early on: $10^{-37}$ seconds 
after the big bang: ``Darn, lets just do it without limits ...!"

\slide{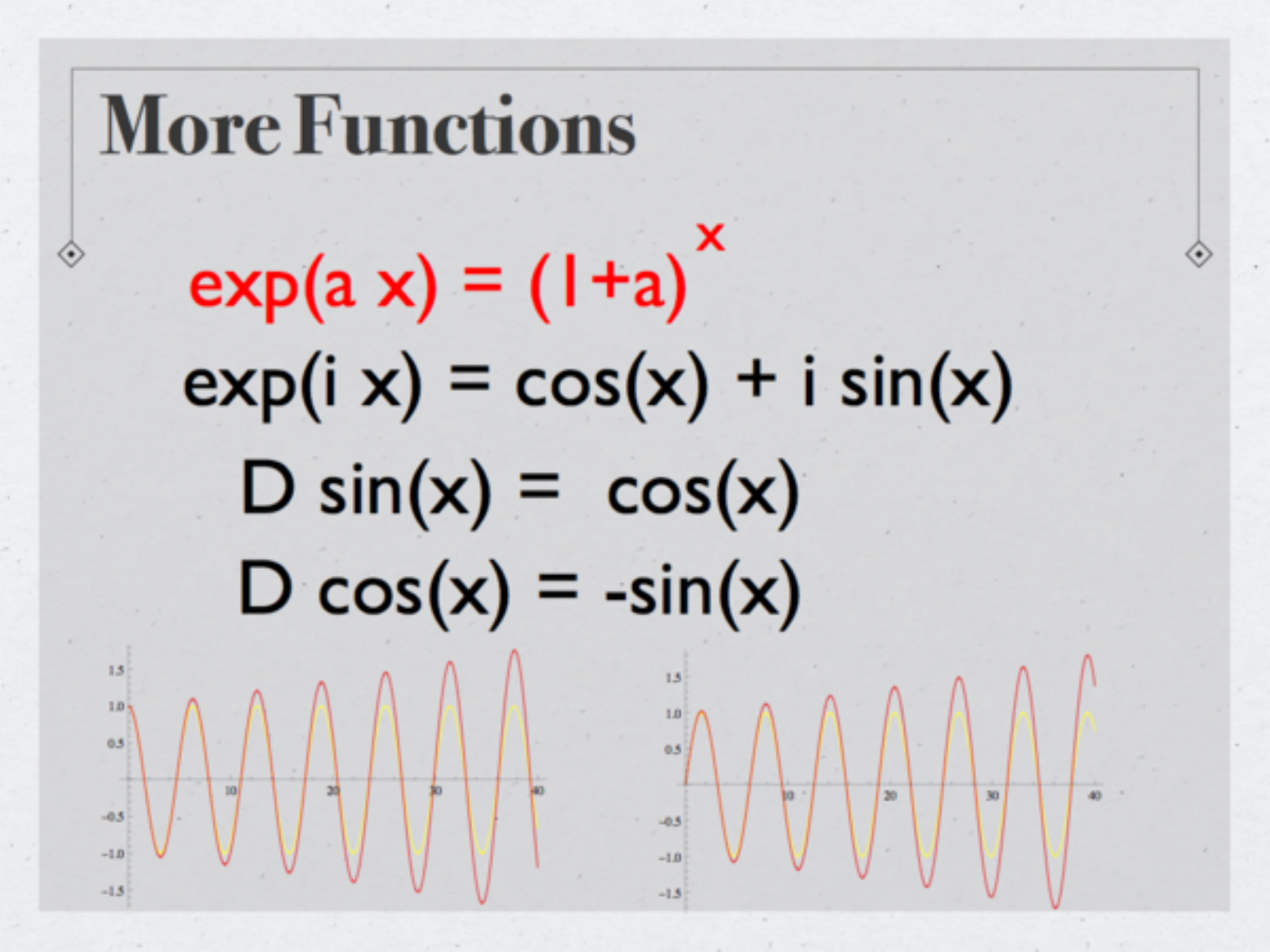}{
We can define cosine and sine by using the exponential, where the interest rate is the square root of -1. 
These deformed trigonometric functions have the property that the discrete derivatives satisfy the
familiar formula as in the books. These functions are close to the trigonometric functions we know if h is small.}

This is by far the most elegant way to introduce {\bf trigonometric functions} also in the continuum.
Unfortunately, due to lack of exposure to complex numbers, it is rarely done. 
The Planck constant is $h=6.626068 \cdot 10^{-34} m^2 kg/s$. 
The new functions $(1+i h a)^{x/h} = \cos(a \cdot x)  + i \sin(a \cdot x)$ are then very close to the 
traditional cos and sin functions. The functions are not periodic
but a growth of amplitude is only seen if $x a$ is of the order $1/h$. 
It needs even for X-rays astronomical travel distances to see the differences. 

\slide{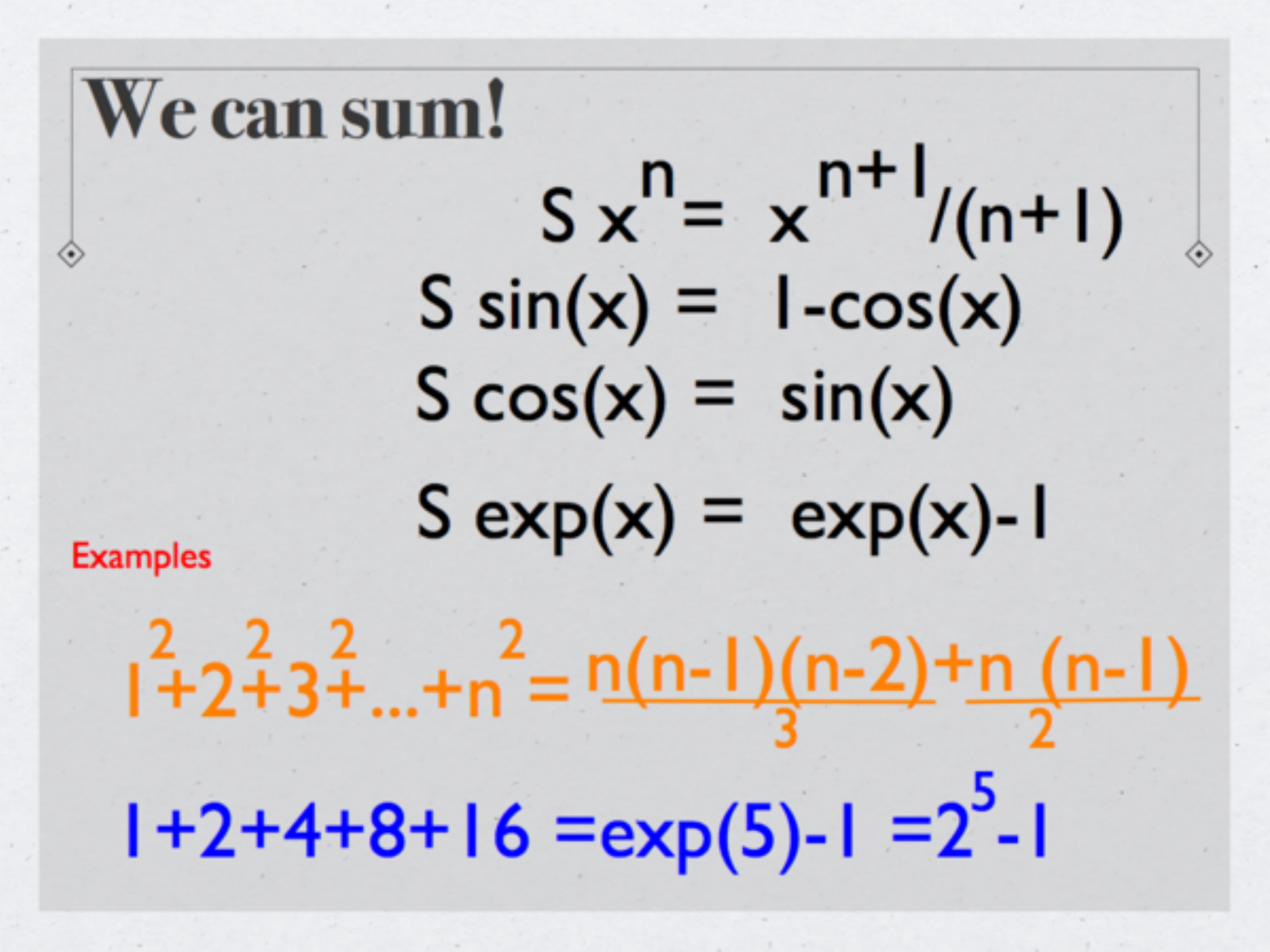}{
The fundamental theorem of calculus is fantastic because it allows us to sum things
which we could not sum before. Here is an example: since we know how to sum up the
deformed polynomials, we can find formulas for the sum of the old squares.
We just have to write the old squares in terms of the new squares which we know how to integrate.}

This leads to a bunch of interesting exercises. For example, because 
$x^2 = [x^2] + [x]$, we have $S x^2 = S [x^2] + S  [x] =  [x^3]/3 +  [x^2]/2$
we get a formula for the sum of the first n-1 squares. Again we have to recall that we sum from 
0 to n-1 and not from 1 to n. 
By the way: this had been a 'back and forth' in the early lesson planning. 
I had started summing from 1 to n and using the backwards
difference $Df(x) = f(x)-f(x-1)$.   The main reason to stick to the forward difference and so to 
polynomials like $[x]^9 = x (x-h) (x-2h)...(x-8h)$
was that we are more familiar with the difference quotient [f(x+h)-f(x)]/h
and also with left Riemann sums. The transition to the traditional calculus becomes easier
because difference quotients and Riemann sums are usually written that way. 

\slide{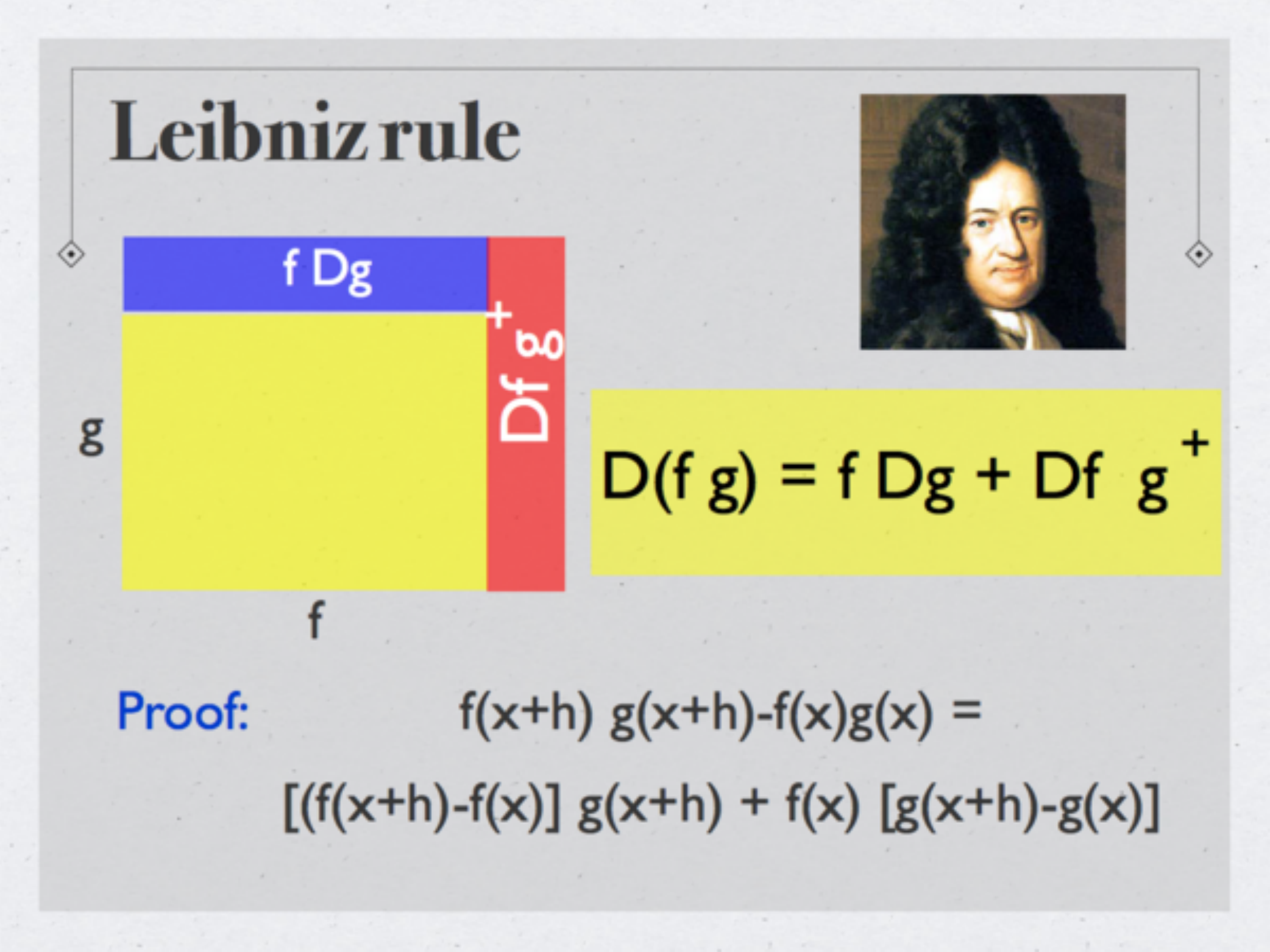}{
We need rules of differentiation and integration. Reversing rules for
differentiation leads to rules of integration. The Leibniz rule is
true for any function, continuous or not. }

The integration rule analogue to integration by parts is called {\bf Abel summation}
which is important when studying Dirichlet series. 
Abel summation is used in this project.
The Leibniz formula has a slight asymmetry. The expansion of the rectangle has two main effects
$f Dg$ and $Df g$, but there is an additional small part $Df Dg$. This is why we have $D(f g) = f Dg + Df g^+$.
The formula becomes more natural when working in the non-commutative algebra mentioned before: if $Df=[s,f]$, then 
$D(fg) = f Dg + Df g$ because the translation operator $s$ now takes care of the additional shift: 
$(Df) g =(f(x+1)-f(x))g(x+1)$. The algebra picture also explains why $[x]^n [x]^m$ is not $[x]^{n+m}$
because also the multiplication operation is deformed in that algebra. 
In the noncommutative algebra we have $(x s^*)^n (x s^*)^m = (x s^*)^{n+m}$. 
While the algebra deformation is natural for mathematicians, it can not be used in calculus courses. 
This can be a reason why it appears strange at first, like quantum mechanics in general. 

\slide{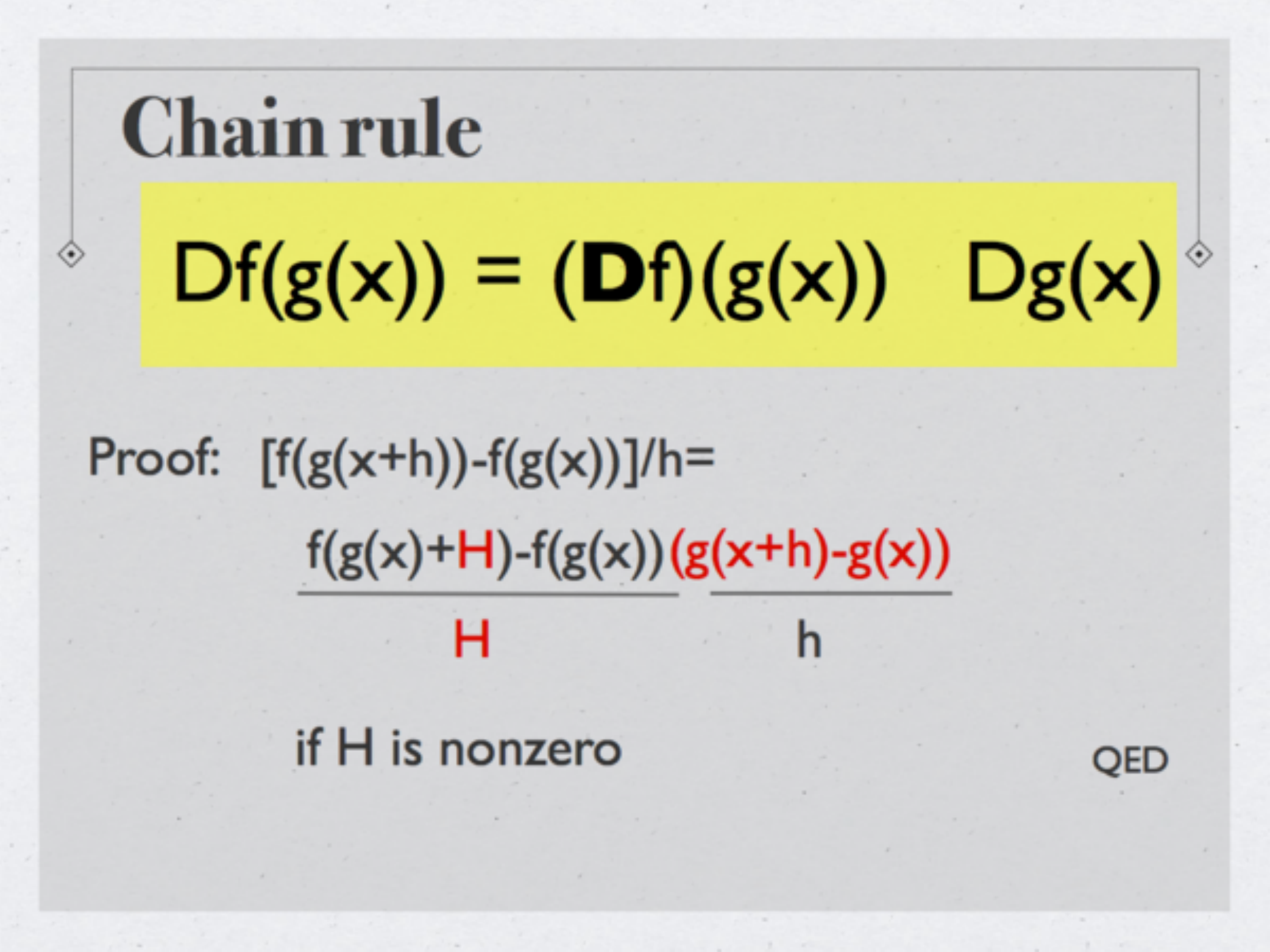}{ 
The chain rule also looks the same. The formula is exact and holds for all
functions. Writing down this formula convinces you that the chain rule is correct. 
The only thing which needs to be done in a traditional course is to take the limit.}

Most calculus textbooks prove the chain rule using linearization: first verify it for linear
functions then argue that the linear part dominates in the limit. 
[Update of December 2013: The proof with the limit \cite{Waestlund}].
Reversing the chain rule leads to the integration tool substitution. 
Substitution is more tricky here because as we see different step sizes h and H in the chain rule. 
This discretization issue looks more
serious than it is. First of all, compositions of functions like $\sin(\sin(x))$ are
remarkably sparse in physics.  Of course, we encounter functions like $\sin(k x)$ but they should be seen 
as fundamental functions and defined like $\sin(k \cdot x) = {\rm Im} (1+k h i)^{x/h}$.
This is by the way different from $\sin(kx) = (1+h i)^{(k x)/h}$. 

\slide{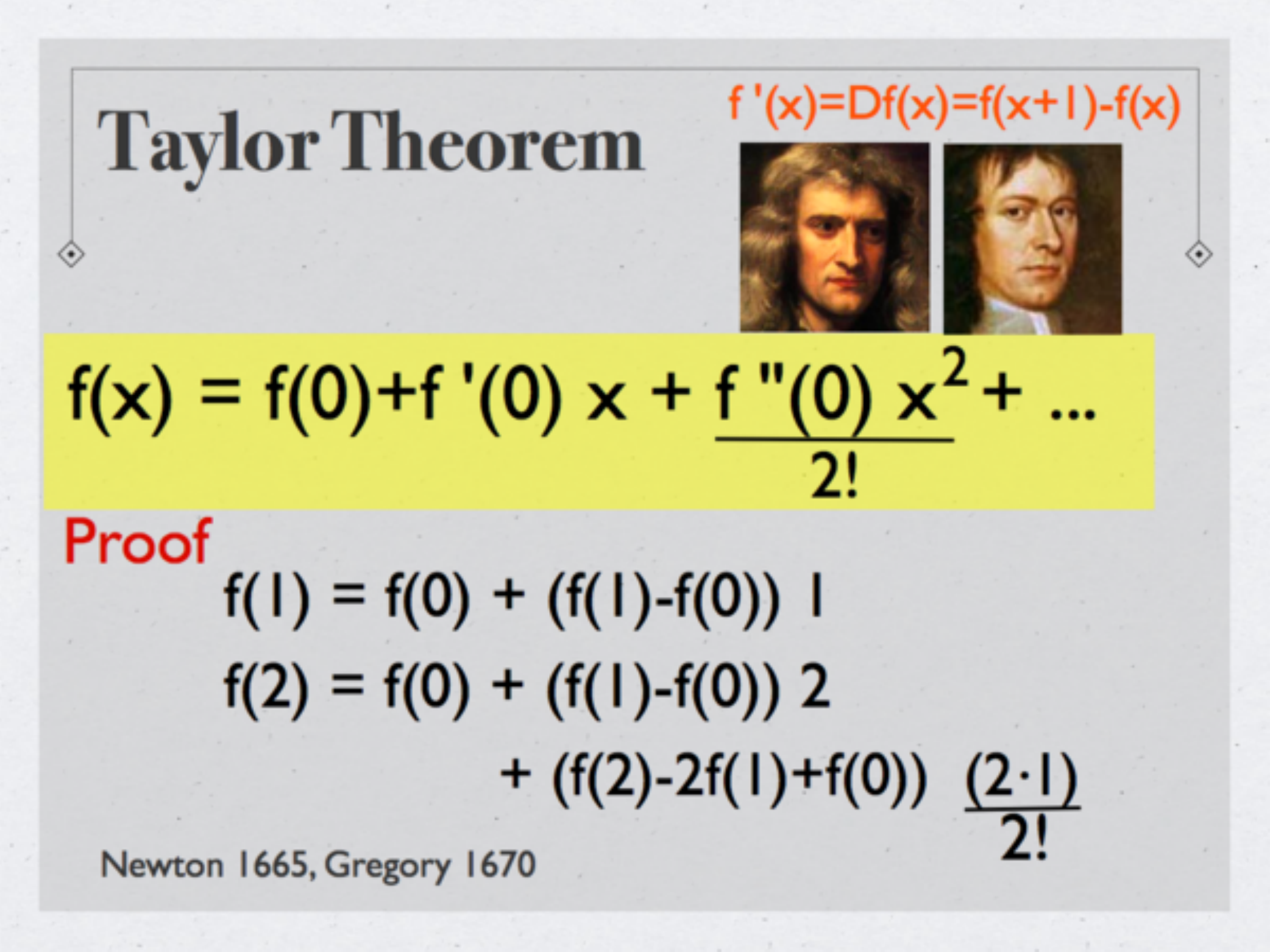}{
Also the Taylor theorem remains true. The formula has been discovered by Newton and Gregory. Written like this,
using the deformed functions, it is the familiar formula we know. We can expand any function, continuous or not, when
knowing the derivatives at a point. It is a fantastic result. You see the start of the proof. It is a nice exercise to 
make the induction step.}

{\bf James Gregory} was a Scottish mathematician who was born in 1638 and died early in 1675.
He has become immortal in the Newton-Gregory interpolation formula.
Despite the fact that he gave the first proof of the fundamental theorem of calculus,
he is largely unknown. As a contemporary of Newton, he was also in the shadow of Newton.
It is interesting how our value system has changed. 
``Concepts" are now considered less important than proofs. 
Gregory would be a star mathematician today. The Taylor theorem can be proven by induction.
As in the continuum, the discrete Taylor theorem can also be proven with PDE methods:
$f(x+t)$ solves the transport equation $D_t f = D f$, where $D_tf(x,t) = f(x,t+1)-f(x,t)$
so that $f(0,t) = [\exp(D t) f](0) = f(0) + Df(0)/1! + D^2f(0)/2! + ...$.

\slide{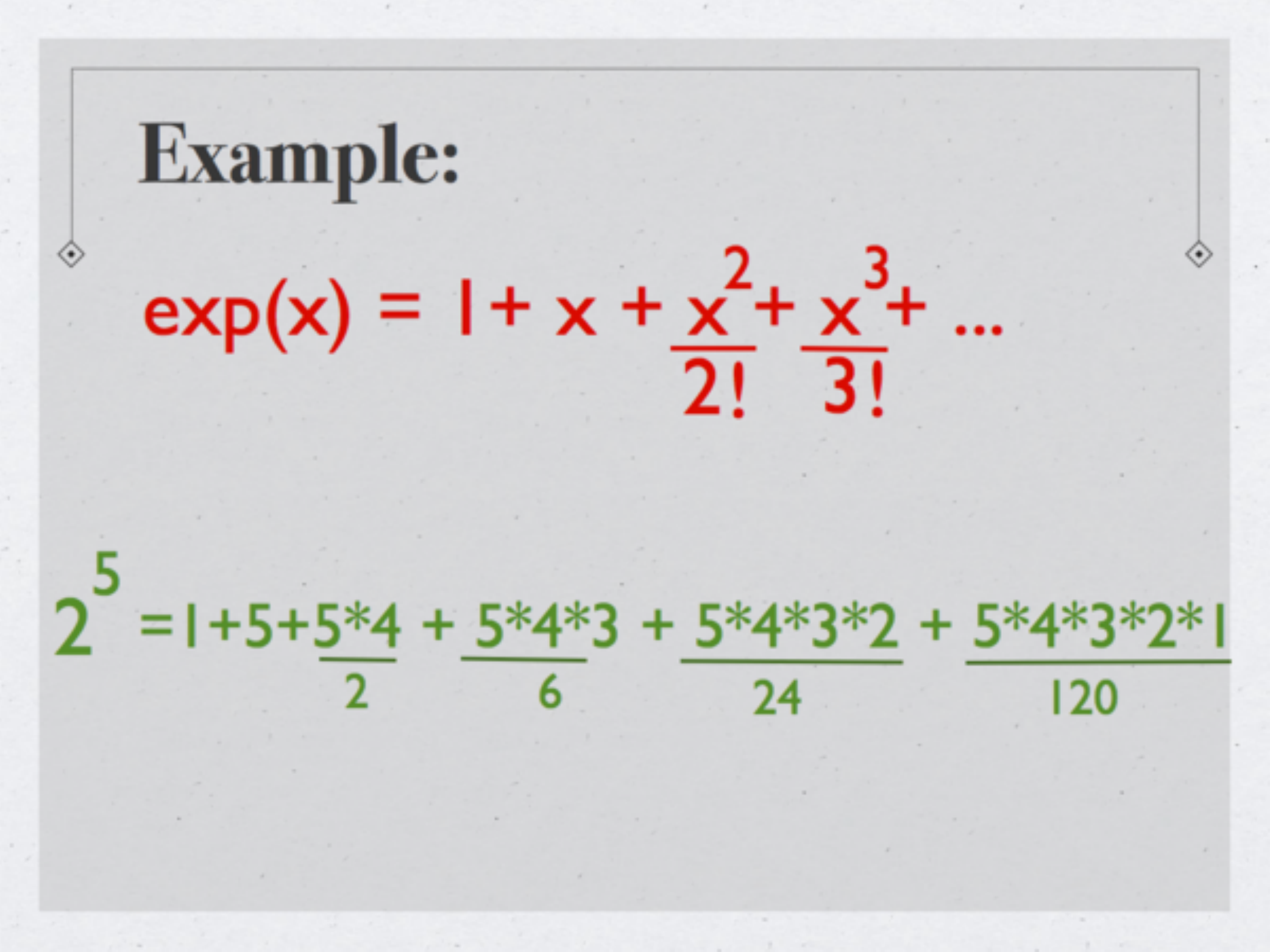}{
Here is an example, where we expand the exponential function and
write it as as a sum of powers. Of course, all the functions are deformed functions. The exponential function
as well as the polynomials $x^n$ were newly defined. We see what the formula means for $x=5$. }

It is interesting that Taylors theorem has such an arithmetic incarnation. 
Usually, in numerical analysis text, this Newton-Gregory result is treated with undeformed 
functions, which looks more complicated.  
The example on this slide is of course just a basic property of the Pascal triangle. 
The identity given in the second line can be rewritten as $32 = 1+5+10+10+5+1$
which in general tells that if we sum the n'th row of the Pascal triangle, we get $2^n$.
Combinatorially, this means that the set of all subsets can be counted by grouping sets of
fixed cardinality. It is amusing to see this as a Taylor formula (or Feynman-Kac). 
But it is more than that: it illustrates again an important point I wanted to make: 
we can not truly appreciate combinatorics if we do not know calculus. 

\slide{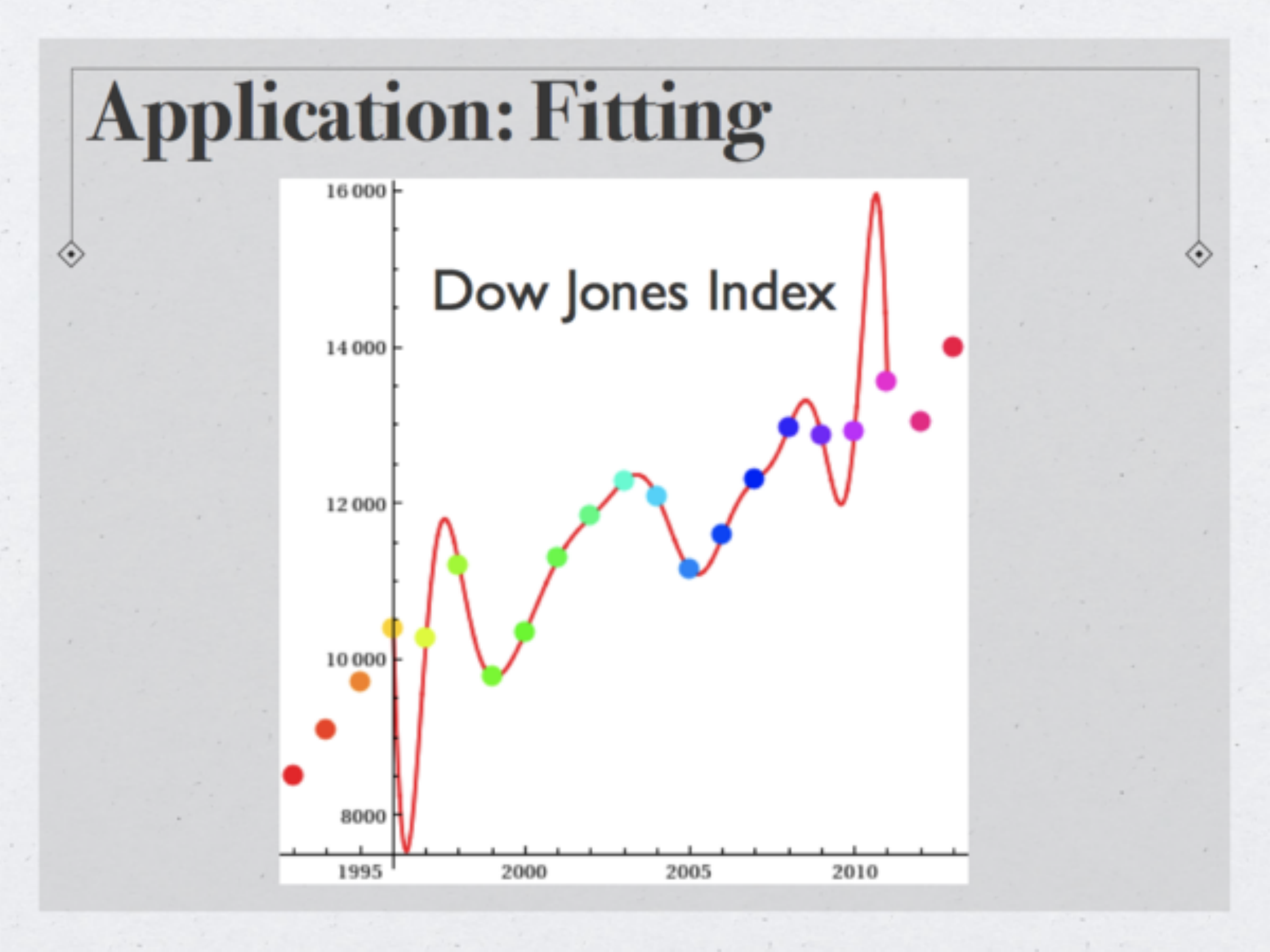}{ 
Taylor's theorem is useful for data interpolation. We can avoid linear algebra and directly write
down a polynomial which fits the data. I fitted here the Dow Jones data of the last
20 years with a Taylor polynomial. This is faster than data fitting using linear algebra.}

When fitting the data, the interpolating function starts to oscillate a lot near the ends. 
But thats not the deficit of the method. If we fit n data points with a polynomial of degree n-1, 
then we get a unique solution. It is the same solution as the Taylor theorem gives.
The proof of the theorem is a simple induction step, using a property of the Pascal triangle. 
Assume, we know $f(x)$. Now apply $D$, to get $Df(x)$. We can now get $f(x+1) = f(x) + Df(x)$. Adding
the terms requires a property of the Pascal triangle. 

\slide{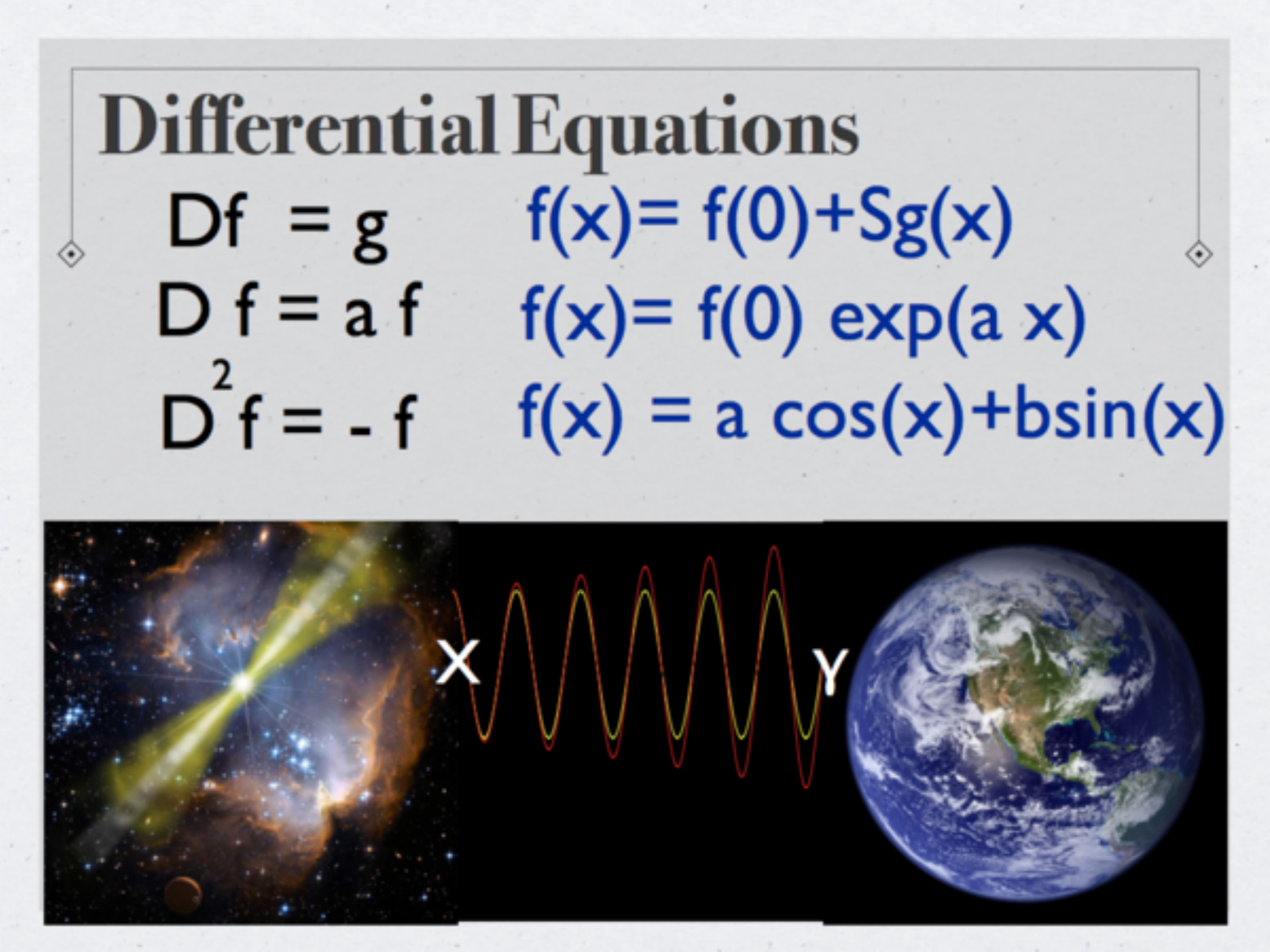}{ 
We can also solve differential equations. We can use the same
formulas as you see in books. We deform the operators and functions so that everything
stays the same. Difference equations can be solved with the same formulas
as differential equations. }

To make the theory stronger, we need also to deform the log, as well as 
rational functions. This is possible in a way so that all the formulas we know in 
classical calculus holds: first define the log as the inverse of exp. Then define
$[x]^{-1} = D \log(x)$ and $[x]^{-n} = D [x]^{1-n}/(1-n)$. 
Now $D x^n = n x^{n-1}$ holds for all integers. 
We can continue like that and define sqrt(x) as the inverse of $x^2$ and 
then the $x^{-1/2}$ as $2D \sqrt{x}$. It is calculus which holds everything together.

\slide{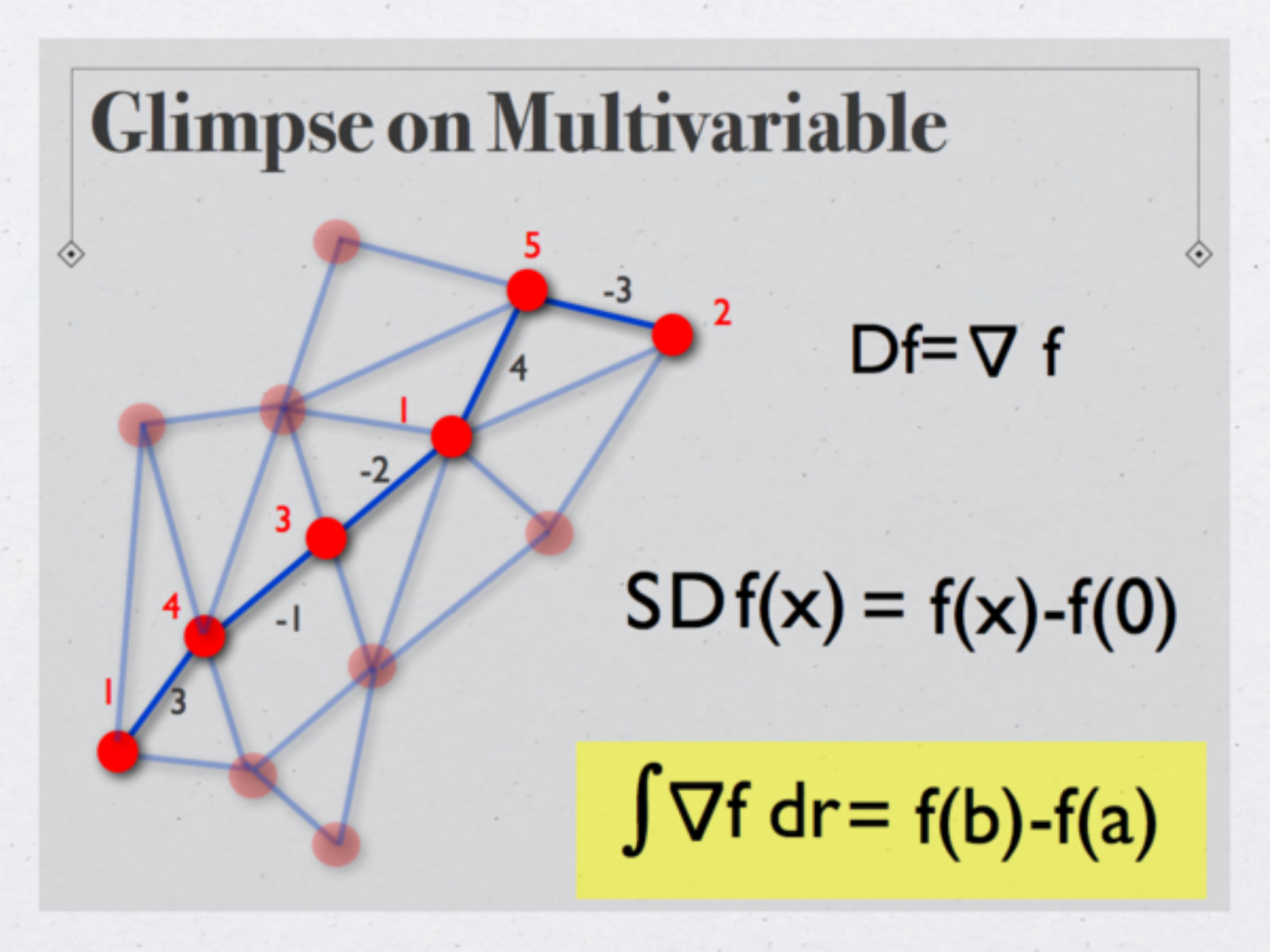}{ 
In multivariable calculus, space takes the form of a graph. Scalar functions are functions on vertices,
Vector fields are signed functions on oriented edges. The gradient of a function is a vector field which is the 
difference of function values. Integrating the gradient of a function along a closed path gives the
difference between the potential values at the end points. This is the fundamental theorem of line integrals. 
}

This result is is also in the continuum the easiest version of Stokes theorem. 
Technically, one should talk about $1$-forms instead of vector fields. 
The one forms are anti-commutative functions on edges. 
\cite{knillcalculus} is an exhibit of three
theorems (Green-Stokes, Gauss-Bonnet \cite{cherngaussbonnet}, 
Poincar\'e-Hopf \cite{poincarehopf}), where everything is defined and 
proven on two pages. See also the overview \cite{classicalstructures}.

\slide{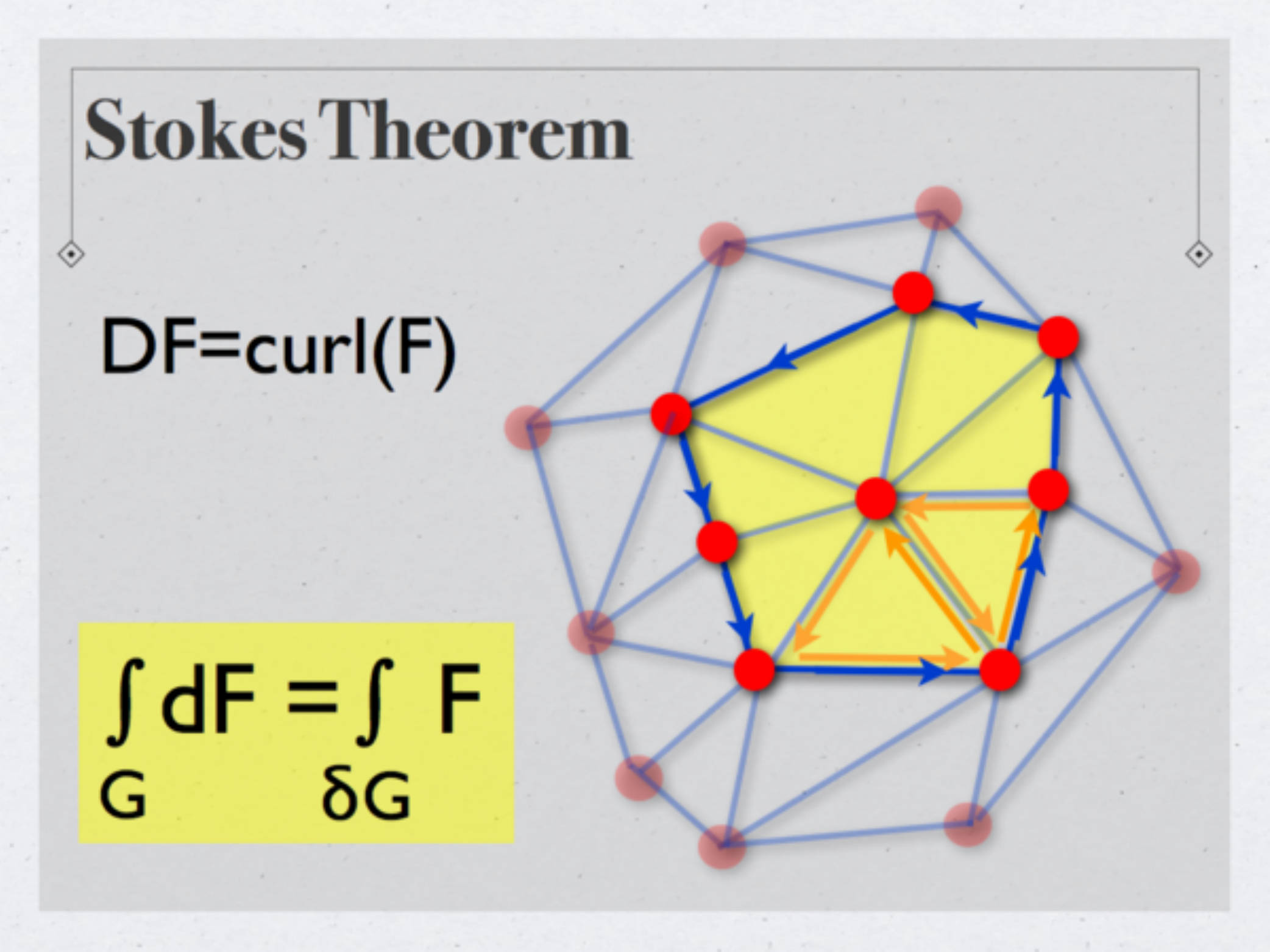}{ 
Stokes theorem holds for a  graph for which the boundary is a graph too.
Here is an example of a "surface", which is a union of triangles. The curl of a vector field F
is a function on triangles is defined as the sum of vector fields along the boundary. Since the terms
on edges in the intersection of triangles cancel, only the line integral along the boundary survives. }

This result is old. It certainly seems have been known to Kirchhoff in 1850. Discrete versions of Stokes pop up again
and again over time. It must have been Poincar\'e who first fully understood the Stokes theorem in all 
dimensions and in the discrete, when developing algebraic topology. He introduced chains because
unlike graphs, chains are closed under boundary operation. This is a major reason, 
algebraic topologists use it even so graphs are more intuitive.
One can for any graph define a discrete notion of differential form as well as an exterior
derivative. The boundary of graph is in general only a chain. For geometric graphs like surfaces made up of
triangles, the boundary is a union of closed paths and Stokes theorem looks like the Stokes theorem we teach. 

\slide{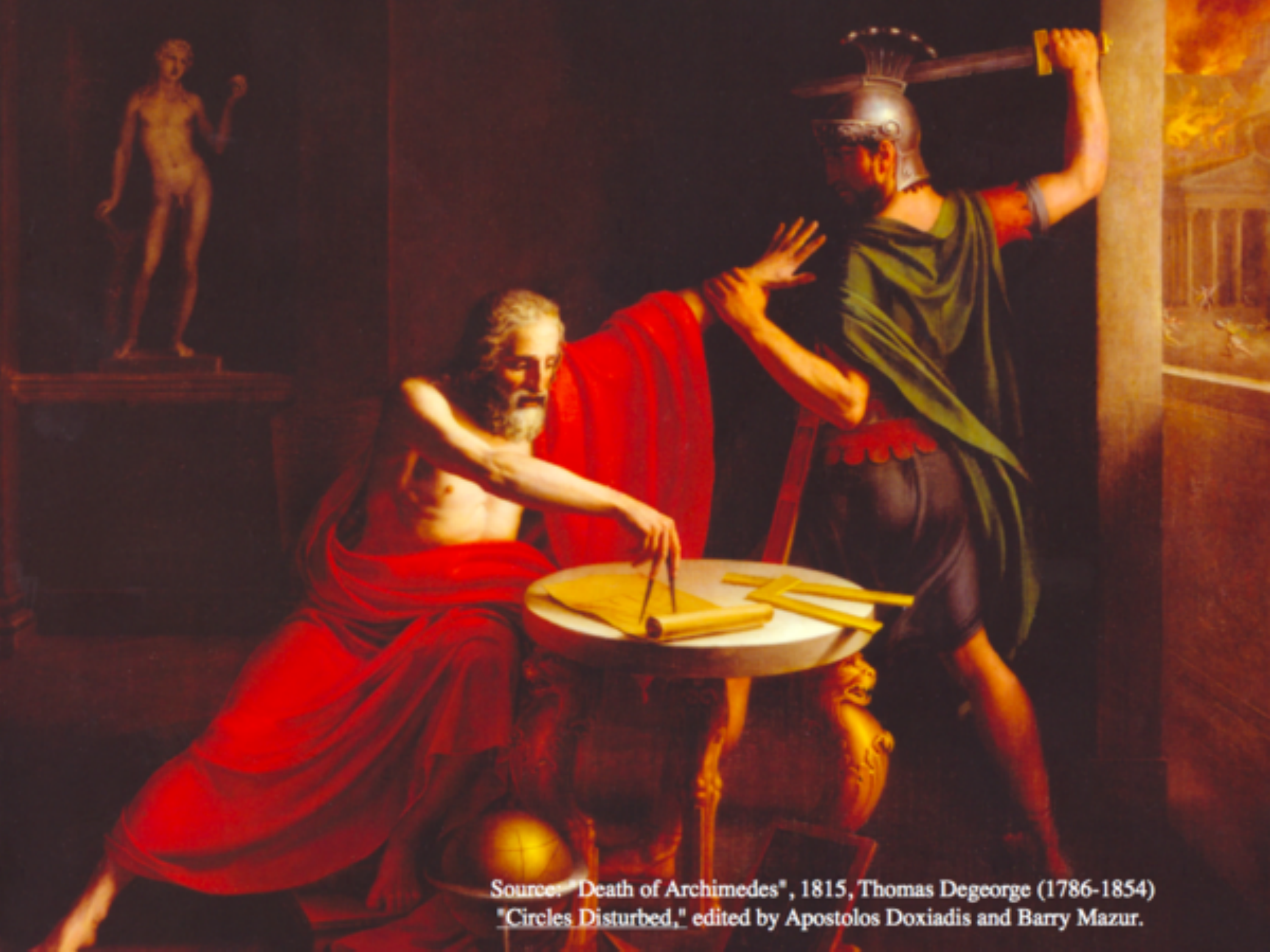}{ 
Could Archimedes have discovered the fundamental theorem of calculus?
His intellectual achievements surpass what we have seen
here by far. \\
Yes, he could have done it. If he would not
have been given the sword, but the concept of a function.}

Both the precise concept
of limit as well as the concept of functions had been missing. While the concept of limit is 
more subtle, the concept of function is easier. The basic ideas of calculus can be explained 
without limits.  Since the ideas of calculus go over so nicely to the discrete, I believe that calculus
is an important subject to teach. It is not only a benchmark and a prototype theory, 
a computer scientist who has a solid understanding of calculus or even 
differential topology can also work much better in the discrete. To walk the talk:
source code of a program in computer vision was written from scratch, takes a movie and finds and tracks
features like corners \cite{KnillHerran}. There is a lot of calculus used inside.  \\

\end{list}
This was a talk given on March 6 at a ``Pecha-Kucha" event at the Harvard Mathematics department, 
organized by {\bf Sarah Koch} and {\bf Curt McMullen}. 
Thanks to both for running that interesting experiment.

\section{Problems in single variable calculus}

Teachers at the time of Archimedes could have assigned problems as follows. 
Some of the following problems have been used in single variable and
extension school classes but most are not field tested. It should become obvious
however that one could build a calculus course sequence on discrete notions.
Since the topic is not in books, it would be perfect also for an ``inquiry based
course", where students develop and find theorems themselves.  \\

\begin{list}{\fcolorbox{color8}{color8}{\Large{\arabic{problem}}}}{\usecounter{problem}}

\problem{
Compute $Dg(n)$ with
$g(n) = [n]^4 = (n-3) (n-2) (n-1) n $
and verify that $Dg(n) = 4 [n]^3$ with $[n]^3=(n-2) (n-1) n$.}

\problem{
Find $Sf(n)$ with $f(n)=n$. This was the task, which the seven year old Gauss was given 
by his school teacher Herr B\"uttner. Gauss computed $Sf(101)=5050$ because he saw that one
can pair 50 numbers summing up to $101$ each and get $[101]^2/2 =101*100/2$. Legend tells that 
he wrote this number onto his tablet and threw it onto the desk of his teacher with the words 
``Hier ligget sie!" (=``Here it is!"). Now its your turn. You are given the 
task to sum up all the squares from $1$ to $1000$. }

\problem{ Find the next term in the sequence
$$2,10,30,68,130,222,350,520,738,1010,1342,... $$
by taking ``derivatives" until you see a pattern, then ``integrate" that pattern. 
}

\problem{ Lets compute some trig function values. We know for example. 
$\sin(0)=0, \sin(1)=1, \sin(2)=2, \sin(3)=1$ and
$\cos(0)=1, \cos(1)=1, \cos(2)=0, \cos(3)=-1$. 
Compute $\sin(10)$ and $\cos(10)$.  ({\bf Answer:} $32$ and $0$.)
}

\problem{
Find $F(n) = Sf(n)$ for the function $f(n) = n^2$. Your explicit formula should satisfy 
$F(1)=1,F(2)=5,F(3)=14$ and lead to the sequence
$$   0,1, 5, 14, 30, 55, 91, 140, 204, 285, \dots \;. $$
}

\problem{
Find $g(x) = Df(x)$ for $f(x) = 3 [x]^5 + 3^x - 2x + 7$ and evaluate it at $x=10$.
{\bf Answer:} We have seen that $D (1+a)^x = a (1+a)^x$ so that $D3^x = 2 \cdot 3^x$. 
$g(x) = 15 [x]^4 + 2 \cdot 3^x - 2$. We have $g(10) = 15 (10*9*8*7) + 2 3^{10} - 2 = 193696$.
Indeed: $f(11) = 343452$ and $f(10) = 193696$ so that the result matches $f(11)-f(10)$. 
}

\problem{
Verify that $Df=f$ for $f=\exp(x)$ and conclude that $Sf(x) = f(x+1)-1$. For
example, for $x=4$ we have $1+2+4+8=16-1$. }

\problem{
The {\bf Fibonnacci sequence} $1,1,2,3,5,8,13,21, \dots$ satisfies the rule $f(x)=f(x-1)+f(x-2)$.
Verify that $Df(x) =f(x-1)$ and use this fact to show that $Sf(x)=f(x+1)-1$. 
For example, for $x=4$ we get $1+1+2+3=8-1$ or for $x=6$ we have $1+1+2+3+5+8=21-1$. 
}

\problem{
Verify that the {\bf Taylor series} $\exp(x) = 1+X+X^2/2!+X^3/3! + \dots$ sums up 
the $n$'th line in the Pascal triangle. For example, $\exp(4)=2^4=1+4+6+4+1$. To do so, check that 
$X^n(x)/n!$ is the {\bf Binomial coefficient} $B(n,x)$ and in particular that $X^n(x)=0$ 
for $n>x$ so that the Taylor series is a finite sum. }

\problem{
Using the definition $\exp(i x)=\cos(x) + i \sin(x)$ of trigonometric functions,
check that $\sin(x) = 1-X^3/3!+X^5/5!-...$ and
$\cos(x) = X-X^2/2!+X^4/4! - ...$. Again, these are finite sums. }

\problem{
{\bf De Moivre} $e^{2 i x} = ({\rm Cos}(x) + i \; {\rm Sin}(x))^2 
= {\rm Cos}(2x) + i \; {\rm Sin}(2x)$ 
is in the continuum the simplest way to derive the 
double angle formulas ${\rm Cos}(2x)={\rm Cos}^2(x) - {\rm Sin}^2(x)$
and ${\rm Sin}(2x) = 2 \; {\rm Cos}(x) {\rm Sin}(x)$. Check whether these double angle
formulas work also with the new functions $\cos(x), \sin(x)$. Be careful that 
$\cos(a \cdot x) = {\rm Re} (1+ia)^x$ is not the same than 
$\cos(ax) = {\rm Re}(1+i)^{ax}$. }

\problem{
Verify that $\sin(x)$ has {\bf roots} at $4 k$ with $k \in Z$ and $\cos(x)$ has roots at
$2+4k$ with $k \in Z$. To do so, note
$\sin(x) = {\rm Im}(e^{x \log(1+i)}) = {\rm Im}(e^{x (\sqrt{2} + i \pi/4)})$. }

\problem{
Find the solution of the {\bf harmonic oscillator} $D^2 f = - 9 f$
which has the initial conditions $f(0)=2,f(1)=5$. Solution: the solution is of the form
$a \cos(3x) + b \sin(3x)$, where $a,b$ are constants. Taking $x=0$ gives $a=2$.
Taking $x=1$ gives $2 \cos(3) + b \sin(3) = 5$. Since $\cos(3)=-2$ and $\sin(3)=2$
we have $b=9/2$. The solution is $f(x) = 2 \cos(3x) + 9 \sin(3x)/2$. }

\problem{
If $h$ is a small positive number called {\bf Planck constant}. Denote by 
${\rm Sin},{\rm Cos}$ the classical $\sin$ and $\cos$ functions. 
Verify that $|\sin(h \cdot (x/h)) - {\rm Sin}(x)| \leq x h$. 
{\bf Solution:} $|{\rm Im} (1+i h)^{x/h} -\Sin(x)| = |{\rm Im} e^{(x/h) {\rm Log}(1+i h)}-{\rm Sin}(x)|$. 
Remark. With $h=7 \cdot 10^{-34}$, a {\bf Gamma ray} of wave length $10^{-12} m$ traveling 
for one million light years $10^{6} \cdot 10^{16}$ meters, we start to see a difference between
$\sin(h \cdot x/h)$ and $\sin(x)$. The amplitude of the former will grow and be seen as
a gamma ray burst. }

\problem{
The functions $\sin(x)$ and $\cos(x)$ are not periodic on the integers.
While we know $\sin(4k)=0$ for every integer $k$, we have $\sin(4k-1)$ exploding
exponentially (this is only since the Planck constant equal to $1$). 
Lets look at the {\bf tangent function} $\tan(x) = \sin(x)/\cos(x)$. Verify that     
it is $4$ periodic $0,1,\infty,-1,0,1,\infty,-1,\dots$ on the integers. }

\problem{
Lets call $x$ be a {\bf discrete gradient critical point} of $f$ if $Df(x)=0$. Assume that $f$ is a continuous
function from the real line $R$ to itself and that has a discrete gradient critical point $a$. Verify that 
there is a classical local maximum or classical local minimum of $f$ in $(a,a+1)$. Solution: this is a 
reformulation of the classical {\bf Bolzano extremal value theorem} on the interval $[a,a+1]$. }

\problem{
Injective functions on $Z$ do not have discrete gradient critical points. Critical points must be defined
therefore differently if one only looks at functions on a graph. 
A point $x$ is a {\bf critical point} if 
the index $i_f(x) = 1-|S(x) \cap \{ f(y) < f(x) \}|$
is not zero. In other words, if either both or no neighbor of $x$ have smaller values than $f(x)$ then $x$ is a
critical point. Verify that if $f$ is an injective $n$ periodic function on $Z$, then $\sum_{k=1}^n i_f(k)$
is equal to $0$. {\bf Hint:} use induction. Assume it is true for $n$, look at $n+1$, then see what happens if
the maximum $x$ of $f$ is taken away. }

\problem{
If $f$ is a classically differentiable function, then for any $a$, there exists
$x \in [a,a+1]$ such that $f'(x) = Df(a)$. Here $f'(x)$ is the classical derivative and
$Df(x)=f(x+1)-f(x)$ the discrete derivative. Solution: this is a reformulation of the 
classical {\bf mean value theorem}. }

\problem{
The {\bf Taylor expansion} $p(x) = \sum_k D^kf(0) X^k/k!$ extends a function $f$ of compact support
from the integers to the real line. It is an explicit {\bf interpolation formula} which provides us
with a function $f(x)$ on the real line satisfying
$p(x)=f(x)$ for integer arguments $x$. Perform the Taylor expansion of the function $f(x)=x^2$.}

\problem{
Check that unlike the classical trigonometric function ${\rm Sin}(x)$ 
the function $\sin(x)$ is not {\bf odd} and the function
$\cos(x)$ is not {\bf even}. Indeed, check that $\sin(x) \to 0$ for 
$x \to -\infty$ while $\sin(x)$ becomes unbounded for $x \to \infty$.
}

\problem{
Compute the sum $\sum_{x=0}^{9} \sin(3 \cdot x)$. 
{\bf Answer:} $-33237  = (1-\cos(3 \cdot 10))/3$. 
Now give an analytic expression for $\sum_{k=0}^{10^{100}} \sin(3 \cdot x)$
which can no more be evaluated by computing the sum directly. }

\problem{
Let us in this problem also write $f'=Df$. 
Integrating the {\bf product rule} $D(f g) = Df g + f^+ Dg$ gives $f g = S (f' g) + S(f^+ g')$
leads to {\bf integration by parts} formula 
$$  S(f' g) = f g - S(f^+ g') \; , $$ 
an expression which is in the discrete also called {\bf Abel summation}. It is usually written as
$$ \sum_{k=0}^{n-1} f'_k g_k = (f_n g_n-f_0 g_0) - \sum_{k=0}^{n-1} f_{k+1} g'_k \; . $$
Use integration by parts to sum up $\sum_{x=0}^{102} \sin(x) x$. 
{\bf Answer:} the sum is $(\sin(103)-\sin(1)) - 102 \cos(103) = -231935380809580545$.
}

\problem{We have defined $1/X$ as the derivative of $\log(x)$. 
Verify that $1/X = {\rm Log}(1+1/x)/{\rm Log}(2)$, where ${\rm Log}$ is
the good old classical natural logarithm. }

\end{list}

\section{Problems in multi variable calculus}

\begin{list}{\fcolorbox{color8}{color8}{\Large{\arabic{problemm}}}}{\usecounter{problemm}}

\problemm{
{\bf Bolzano's extremal value theorem} is obvious for finite graphs because a function on a finite set 
always has  a global maximum and minimum. Verify that for a global minimum, the 
{\bf index} $i_f(x)$ is equal to $1$. For every integer $n$ find a graph $G$ and a function $f$ for which 
a maximum is equal to $n$. For example, for a star graph with $n$ rays and a function $f$ having the
maximum in the center $x$, we have $i_f(x)=1-n$. }

\problemm{
A {\bf curve} in $G$ is a subgraph $H$ with the property that every vertex has either $1$
or $2$ neighbors. The {\bf arc length} of a path is the number of edges in $H$. The boundary
points are the vertices with one neighbor.
A connected {\bf surface} $H$ in $G$ is a subgraph for which every vertex has either 
a linear graph or circular graph as a unit sphere. The surface area of $H$ is the number of 
triangles in $H$. The interior points are the points for which the unit sphere is a circular graph.
A surface is called {\bf flat} if every interior point has a unit sphere which is $C_6$. Flat surfaces play
the role of regions in the plane. The curvature of a boundary point of a flat surface is $1-|S(x)|/6$. 
Verify that the sum of the curvatures of boundary points is $1-g$ where $g$ is the number of holes. 
This is a discrete {\bf Hopf Umlaufsatz}. 
}

\problemm{
A graph $G$ is called {\bf simply connected} if every simple and connected curve $C$ in $G$ can be deformed to a single 
edge using the following two {\bf deformation steps} of a $C$: replace two edges of $C$ within a triangle with the 
third edge of that triangle or replace an edge of $C$ in a triangle with the two other edges of the triangle. 
A deformation step is called {\bf valid} if the curve remains a simple connected curve or degenerates to $K_3,K_2$
after the deformation. A $1$-form $F$ is called a {\bf gradient field}, if $F=df$ for some scalar function $f$. 
Verify that on a simply connected surface, a $1$-form $F$ is a gradient field if and only if the curl $dF$
is zero on each triangle. }

\problemm{
{\bf Stokes theorem} for a flat orientable surface $G$ is what traditionally is 
called {\bf Greens theorem}. It tells that the sum of the curls over the triangles 
is equal to the line integral along the boundary. Look at the wheel graph $W_6$ which
is obtained by adding an other vertex $x$ to a circular graph $C_6$ and connecting all 
vertices $y_i$ of $C_6$ with $x$. Define the $1$-form $F$, where each of the central edges $(x,y_i)$
is given the value $1$ and each edge $(y_i,y_{i+1})$ is given the value $i$. Find the curl $dF$
on all the triangles and verify Green's theorem in that case. }

\problemm{
A graph $G$ is a {\bf surface} if for every vertex $x \in G$, the unit sphere $S(x)$ is a circular graph. 
The curvature at a point $x$ of a surface is defined as $1-|S(x)|/6$. Verify that an octahedron is a
discrete surface and that the curvature is constant $1/3$ everywhere. See that the curvatures add up to $2$. 
Verify also that the icosahedron is a discrete surface and that the curvature is constant $1/6$
adding up to $2$. You have verified two special cases of {\bf Gauss-Bonnet}. Play with other discrete surfaces
but make sure that all ``faces" are triangles. 
}

\problemm{
We want to see an analogue of the {\bf second derivative test} for functions $f(x,y)$ of two variables
on a discrete surface $G$ like an icosahedron. For an injective function $f$ on the vertices, call $x$ a 
{\bf critical point} with nonzero index if $S^-(x) = S(x) \cap \{ y \; | \; f(y)<f(x) \}$ has 
Euler characteristic different from $1$. The number $i_f(x)=1-\chi(S^-(x))$ plays the role of the {\bf discriminant}. 
Verify that if it is negative then we deal with a {\bf saddle point} which is neither maximum nor minimum. 
Verify that if it is positive with empty $S^-(x)$, then it is a {\bf local minimum}. Verify that if 
the discriminant is positive and $S^-(x)$ is a circular graph, then it is a {\bf local maximum}. 
As in the continuum, where the discriminant at critical points can be zero, this can happen also here but
only in higher dimensions. Maybe you can find an example of a function on a three dimensional graph. }

\problemm{
Given a surface $G=(V,E)$ and an injective function $f$ on $V$,
define the {\bf level curve} $C$ of $f=c$ as the graph which has as vertices 
the edges of $G$ for which $f$ changes. 
Two such vertices are connected if they are contained in 
a common triangle. (Graph theorists call this a subgraph of the {\bf ``line graph"} of $G$. 
Verify that for a two dimensional geometric graph,  level curve $C$ of $f$ 
is a finite union of circular graphs and so a geometric one dimensional graph. 
The fact that the gradient $df(x)$ is a function on $C$ changing sign 
on every original edge is the analogue of the classical fact that the gradient is perpendicular
to the level curve. }

\problemm{
An oriented edge $e$ of a graph $G$ is a "direction" and plays the role of unit vectors
in classical calculus. The value of $df(e)$ is the analogue of a {\bf directional 
derivative}. Start with a vertex $x$ of $G$, then go into direction $e$, 
where the directional derivative $df(e)$ is maximal and stop if there 
is no direction in which the directional derivative is positive. 
This discrete {\bf gradient flow} on the vertex set leads to a local maximum. 
Investigate how this notion of directional derivative corresponds to the notion 
$D_vf(x) = \nabla f(x) \cdot v$ in the continuum. Especially note that for an 
injective function, the directional derivative is never zero.  }

\problemm{
An injective function $f$ on a graph has no points for which $df(a)=0$. We therefore have 
to define critical points differently. We say $x$ is a {\bf critical point} of $f$ 
if $i_f(x) = 1-\chi(S_f^-(x) \})$ is not zero. Let $G$ be a discretisation of a disc. 
Check that for maxima, $i_f(x)=1$ and $S_f^-(x)$
is a circular graph. Check that for minima, $i_f(x)=1$ and $S_f^-(x)$ is empty. For
Saddle points, $i_f(x)=-1$ where $S_f^-(x)$ consists of two linear graphs. 
A vertex is called a {\bf monkey saddle} if $i_f(x)=-2$. 
Verify the "island theorem": the sum over all indices over the disc is equal to $1$.
A discrete {\bf Poincar\'e-Hopf} theorem \cite{poincarehopf} assures that in general 
$\sum_{x} i_f(x) = \chi(G)$.
}

\problemm{
Let $F,G$ be two $1$-forms, define the {\bf dot product} of $F,G$ at a vertex $x$ 
as $\langle F,G \rangle(x) = \sum_{e \in E(x)} F(e) G(e)$, where $E(x) \subset E$ is 
the set of edges attached to $x$.  Define the length $|F|(x)$ as 
$\sqrt{\langle F,F \rangle}$. The {\bf angle} between two $1$-forms at $x$ is defined
by the identity $\langle F,G \rangle = |F| \cdot |G| \cos(\alpha)$. 
Given two functions $f,g$ and an edge $e=(a,b)$ define the angle between two
{\bf hyper surfaces} $f=c,g=d$ at an edge $e$ as the angle between $df(e)$ and $dg(e)$.
The hyper surface $f=c$ is part of the {\bf line graph} of $G$: 
it has as vertices the edges of the graph at which $f$ changes sign
and as edges the pairs of new vertices which intersect when seen as old edges. }

\problemm{
Let $F,G$ be two $1$-forms. Define the {\bf cross product} of $F,G$ at a vertex
$x$ as the function on adjacent triangle $(x,y,z)$ as 
$F \times G(x,y,z) = F(x,y) G(x,z) - F(x,z) G(x,y)$. 
It is a number which is the analogue of projecting the classical cross product onto a 
$2$-dimensional direction in space. The cross product is anti-commutative. Unlike a two form, it is
not a function on triangles because the order of the product plays a role. It can be made
a $2$-form by anti-symmetrization. Here is the problem: verify that for any two functions 
$f,g$ on the vertices of a triangle the cross product $df \times dg$
is up to a sign at all vertices. }

\problemm{
Given three $1$-forms $F,G,H$, define the {\bf triple scalar product} at a vertex $x$
as a function on an adjacent tetrahedron $(x,y,z,w)$ as
$F \cdot G \times H(x,y,z,w) = {\rm det}( \left[ \begin{array}{ccc} 
                                    F(x,y)&F(x,z)&F(x,w) \\ 
                                    G(x,y)&G(x,z)&G(x,w) \\
                                    H(x,y)&H(x,z)&H(x,w) \end{array} \right] )$. 
Verify that this triple scalar product has the same properties as the triple scalar product 
in classical multivariable calculus. }

\problemm{
Let $G=K_2$ be the complete graph with two vertices and one edge. Verify that the Dirac operator $D$
and form Laplacian $L=D^2$ are given by
$$  D=\left[ \begin{array}{ccc} 0 & 0 & -1 \\ 0 & 0 & 1 \\ -1 & 1 & 0 \\ \end{array} \right],
    L=\left[ \begin{array}{ccc} 1 & -1 & 0 \\ -1 & 1 & 0 \\ 0 & 0 & 2 \\ \end{array} \right] \; . $$
For an informal introduction to the Dirac operator on graphs and Mathematica code, see \cite{DiracKnill}. 
The Laplacian has the eigenvalues $2,2,0$. Write down the solution of the wave equation $f''=-L f$
where $f(t) = (u(t),v(t),w(t))$ and $u(t),v(t)$ are the values of $f$ on vertices and $w(t)$
is the evolution of the $1$-form. 
}

\problemm{
The graph G=$K_5$ is the smallest possible four dimensional space. 
Given a $1$-form $A$, we have at each vertex $6$ triangles. The $2$-form 
$F=dA$ defines so a field with $6$ components. If one of the four edges
emerging from $x$ is called "time" and the $3$ triangles containing this edge
define the electric field $E_1,E_2,E_3$ and the other $3$ triangles are
the magnetic field $B_1,B_2,B_3$. The equations $dF=0,d^*F=j$ are
called {\bf Maxwell equations}. Can we always get $F$ from the current $j$?
Verify that if ${\rm div}(A) = d^* A=0$ ({\bf Coulomb gauge}), 
then the {\bf Poisson equation} $LA=j$ always allows to compute $A$ and so $F=dA$.
The Coulomb gauge assumption is one of Kirchhoff's conservation law: the sum 
over all in and outgoing currents at each node are zero. }

\problemm{
Let $G=C_5$ denote the circular graph with $5$ vertices. Find a function on 
the vertex set, which has exactly two critical point, one minimum of index
$1$ and one maximum of index $-1$. Now find a function with $4$ critical points
and verify that there are no functions with an odd number of critical points. 
{\bf Hint:} show first that the index of a critical point is either $1$ or $-1$ and
use then that the indices have to add up to $0$. 
}

\problemm{
Now look at the tetrahedron $G=K_5$. Verify that every function on $G$
has exactly one critical point, the minimum. 
}

\problemm{
Take a linear graph $L_3$ which consists of three vertices and two edges.
There are 6 functions which take values $\{1,2,3 \}$ on the vertices. 
In each case compute the index $i_f(x)$ and verify that the average 
over all these indices gives the curvature of the graph, which is $1/2$
at the boundary and $0$ in the center. This is a special case of a result
proven in \cite{indexexpectation}.  }

\end{list}

\section*{Appendix A: Graph theory glossary}

A {\bf graph} consists of a vertex set $V$ and an edge set $E \subset V \times V$. 
If $V$ and $E$ are finite, it is called a {\bf finite graph}. If every $e=(a,b) \in E$
satisfies $a \neq b$, elements $(a,b)$ and $(b,a)$ are identified and every $e$ appears
only once, the graph is called a {\bf finite simple graph}. In other words, we 
{\bf disregard loops} $(a,a)$, look at {\bf undirected graphs} without 
{\bf multiple connections}. The cardinality of $V$ is called {\bf order} and the 
cardinality of $E$ is called {\bf size} of the graph. When counting $E$, we look at the
unordered pairs $e=\{a,b\}$ and not the ordered pairs $(a,b)$. The later would be twice
as large.  \\

\begin{tiny}
Restricting to finite simple graphs is not much of a restriction
of generality. The notion of loops and multiple connections can be incorporated later by
introducing the algebraic notion of {\bf chains} defined over the graph, directions using orientation.
Additional structure can come from divisors or difference operators defined over the graph. 
When looking at geometry and calculus in particular, the notion of chains resembles more the
notion of fiber bundles, and directions are more like a choice of {\bf orientation}. 
Both enrichments of the geometry are important but can be popped upon a finite simple graph if
needed, similarly as bundle structures can be added to a manifold or variety. 
As in differential geometry or algebraic geometry, the notions of fiber bundles or 
divisors are important because they allow the use of algebra to study geometry. Examples are
cohomology classes of differential forms or divisor classes and both work very well in graph theory too.  \\
\end{tiny}

Given a finite simple graph $G=(V,E)$, the elements in $V$ are called {\bf vertices} or {\bf nodes}
while the elements in $E$ are called {\bf edges} or {\bf links}. A finite simple graph contains more
structure without further input. A graph $H=(W,F)$ is called a {\bf subgraph} of $G$ if
$W \subset V$ and $F \subset E$. A graph $(V,E)$ is called {\bf complete} if the size is $n(n-1)/2$,
where $n$ is the order. In other words, in a complete graph, all vertices are connected. The
{\bf dimension} of a complete graph is $n-1$ if $n$ is the order. A single vertex is a complete graph
of dimension $0$, an edge is a complete graph of dimension $1$ and a {\bf triangle} is a complete graph
of dimension $2$. Let $\G_k$ denote the set of complete subgraphs of $G$ of dimension $k$. 
We have $\G_0 = V, \G_1 = E$ and $\G_2$ the set of triangles in $G$ and $\G_3$ the set of tetrahedra in $G$. 
Let $v_k$ be the number of elements in $\G_k$. 
The number $\chi(G) = \sum_{k=0}^{\infty} (-1)^k v_k$ is called the 
{\bf Euler characteristic} of $G$. For a graph without tetrahedra for example, we have
$\chi(G) = v-e+f$, where we wrote $v=v_0, e=v_1, f=v_2$. For example, for an octahedron or
icosahedron, we have $\chi(G)=2$. More generally, any triangularization of a two dimensional sphere
has $\chi(G)=2$. The cube graph has no triangles so that $\chi(G) = 8-12=-4$ illustrating that we see
it as a sphere with $6=2g$ holes punched into. For a discretisation of a sphere with $2g$ holes, we 
have $\chi(G) = 2-2g$. If $g=1$ for example, where we have $2$ holes, we can identify the boundaries
of the holes without changing the Euler characteristic and get a torus.  \\

\begin{tiny}
Often, graphs are considered as "curves", one dimensional objects in which higher dimensional 
structures are neglected. The most common point of view is to disregard even the two-dimensional 
structures obtained from triangles and define the Euler characteristic of a graph as 
$\chi(G) = v-e$, where $v$ is the number of vertices and $e$ the number of edges. 
For a connected graph this is then often written as $1-g$, where $g$ is called the {\bf genus}.
It is a linear algebra computation to show that the Euler characteristic $\sum_k (-1)^k v_k$
is equal to the cohomological Euler characteristic $\sum_k (-1)^k b_k$, where $b_k$ is the
dimension of the cohomology group $H^k(G)$. For example, $b_0$ is the number of connected
components of the graph and $b_1$ is the `number of holes" the {\bf genus} of the curve. 
Ignoring higher dimensional structure also means to ignore higher dimensional cohomologies
which means that the Euler-Poincar\'e formula is $v_0-v_1 = b_0 - b_1$ or $v-e = 1-g$. \\
\end{tiny}

Given a subset $W$ of the vertex set, it {\bf generates} a graph $H=(W,F)$ which has
as edge set $F$ all the pairs $(a,b)$ of $W$ with $(a,b) \in E$. A vertex $y$ is called
a {\bf neighbor} of $x$ if $(x,y) \in E$. 
If $x$ is a vertex of $G$, denote by $S(x)$ the {\bf unit sphere} of of $x$. It is
the graph generated by the set of neighbors of $x$.  Given a subgraph $H=(W,F)$ 
of $G=(V,E)$, we can build a new graph $G=(V \cup \{v\},E \cup \{ (v,y) \; | \; y \in W \})$
called {\bf pyramid extension} over $H$. We can define a class of {\bf contractible graphs}
as follows: the one point graph is contractible. Having defined contractibility for graphs
of order $n$, a graph $G$ of order $n+1$ is {\bf contractible} if there exists a contractible subgraph
$H$ of $G$ such that $G$ is a pyramid extension over $H$. 
Having a notion of contractibility allows to define {\bf critical point} properly.
And critical points are one of the most important notions in calculus. 
Given an injective function $f$ on the vertex set $V$, a vertex $x$ is called a 
{\bf critical point} of $f$ if $S^-(x) = S(x) \cap \{ y \; | \; f(y)<f(x) \}$ is
empty or not contractible. Because contractible graphs have Euler characteristic $1$,
vertices for which $i_f(x)$ is nonzero are critical points. As in the continuum, there
are also critical points with zero index.  \\

Lets look at examples of graphs: the {\bf complete graph} $K_{n+1}$ is an $n$ dimensional
simplex. They all have Euler characteristic $1$. 
The {\bf circular graph} $C_n$ is a one-dimensional geometric graph. The Euler characteristic is $0$. 
Adding a central vertex to the circular graph and connecting them to all the $n$ vertices in $C_n$
produces the {\bf wheel graph} $W_n$. If the original edges at the boundary are removed, we get the
{\bf star graph} $S_n$. All the $W_n$ and $S_n$ are contractible and have Euler characteristic $1$. 
The {\bf linear graph} $L_n$ consists of $n+1$ vertices connected with $n$ edges. Like the star graph it
is also an example of a tree. A {\bf tree} is a graph graph without triangles and subgraphs $C_n$. 
Trees are contractible and have Euler characteristic $1$. The {\bf octahedron} and {\bf icosahedron} 
are examples of two dimensional geometric graphs. Like the sphere, they have Euler characteristic $2$. 

\section*{Appendix B: Historical notes}

For more on the history of calculus and functions and people, 
see \cite{Cajori,Bell,BoyerCalculus,EdwardsHistory,Bottazzini,DunhamMonthly,KatzHistoria,Arnold,Rosenthal,Bardi,Tent,Dunham}. \\

{\bf Archimedes} (287-212 BC) is often considered the father of {\bf integral calculus} as he introduced
Riemann sums with equal spacing. Using comparison and
exhaustion methods developed earlier by {\bf Eudoxos} (408-347 BC), he 
was able to compute volumes of concrete objects
like the hoof, domes, spheres or areas of the circle or regions below
a parabola. Even so the concept of function is now very natural to us,
it took a long time until it entered the mathematical vocabulary. 
\cite{Bottazzini} argues that tabulated function values used in ancient astronomy
as well as tables of cube roots and the like can be seen as evidence that
some idea of function was present already in antiquity.  It was
{\bf Francois Vi\`ete} (1540-1603) who introduced variables leading to
what is now called ``elementary algebra" and
{\bf Ren\'e Descartes} (1596-1650) who introduced analytical geometry as well as
{\bf Pierre de Fermat} (1601-1665) who studied also extremal problems and
{\bf Johannes Kepler} (1571-1630) who acceleration as rate of change of
velocity and also optimized volumes of wine casks or computed what we would call today
elliptic integrals or trigonometric functions. \\

As pointed out for example in \cite{DieudonneMusic}, new ideas are
rarely without predecessors. The notion of coordinate for example has been anticipated
already by {\bf Nicole Oresme} (1323-1392) \cite{BoyerCalculus}. 
Dieudonn\'e, who praises the idea of function as a ``great
advance in the seventeenth century", points out that first examples of plane curves
defined by equations can be traced back to the fourth century BC. 
While ancient Greek mathematicians already had a notion 
of tangents, the notion of derivative came only with {\bf Isaac Barrow} (1630-1703) 
who computed tangents with a method we would today called implicit differentiation.
{\bf Isaac Newton} introduced derivatives on physical grounds using different terminology 
like ``fluxions". Until the 18th century, functions were treated as an analytic expression
and only with the introduction of set theory due to {\bf Georg Cantor}, the modern 
function concept as a single-valued map from a set $X$ to a set $Y$ appeared. 
From \cite{Luzin}:
"The function concept is one of the most fundamental concepts of modern
mathematics. It did not arise suddenly. It arose more than two hundred years ago
out of the famous debate on the vibrating string".
The word ``function" as well as modern calculus notation came to us through 
{\bf Gottfried Leibniz} (1646-1716) in 1673. (\cite{Bottazzini} tells that
the term ``function" appears the first time in print in 1692 and 1694).
Trig functions entered calculus with Newton in 1669 but were included in textbooks only in 
1748 \cite{KatzHistoria}. {\bf Peter Dirichlet} (1805-1859) gave
the definition still used today stressing that for the definition $x \to f(x)$ 
it is irrelevant on how the correspondence is established. 
The notation $f(x)$ as a value of $f$ evaluated at a point $x$ came even later through
{\bf Leonard Euler} (1707-1783) and functions of several variables only appeared
at the beginning of the eighteenth century.  \\

When Fourier series became available with {\bf Jean Baptiste Fourier} (1768-1830) \cite{FourierChaleur},
one started to distinguish between ``function" and ``analytic expression" for
the function. The concept has evolved ever since. One can see a function as a subset $Z$ 
of the Cartesian product $X \times Y$ such that $(\{x\} \times Y) \cap Z$ has exactly one point. 
A more algorithmic point of view emerged with the appearance of
computer science pioneered by Fermat and {\bf Blaise Pascal} (1623-1662).
It is the concept of a function as a ``rule" assigning to an ``input" $x$ an ``output" $y$. 
A notion of Turing machine made precise what it means to ``compute" $f(x)$. 
After the notion of {\bf analytic continuation} has sunk in by 
{\bf Bernhard Riemann} (1826-1866), it became clear that a function can be a world by itself 
with more information than anticipated. The zeta function $\zeta(z)$ for example,
given as a sum $\zeta(s) = \sum_{n=1}^{\infty} n^{-s}$ can be extended beyond
the domain $Re(s)>1$, where the sum is defined to the entire complex plane except $s=1$. 
We already know in school calculus that the {\bf geometric sum} $f(a) =\sum_{n=0}^{\infty} a^n$ 
makes sense for $a>1$ if written as $f(a) = 1/(1-a)$ defining so $f(2)=-1$ even if
the sum $\sum_n 2^n$ diverges. The theory of Riemann surfaces shows that functions often 
have an extended life as geometric objects. We have also learned that there are functions which
are not computable. The function which assigns to a Turing machine the value $1$ if it halts and 
$0$ if not for example is not computable.  The geometric point of view is how everything has started:
functions were perceived as geometry very early on. The relation $f(x,y,z)=x^2+y^2-z^2=0$ for example
would have been seen by Pythagoras as a relation of areas of three squares built from a right angle triangles 
and mathematicians like {\bf Brahmagupta} (598-670) or {\bf Al-Khwarizmi} (780-850) 
saw the construction of roots of the quadratic equation $f(x) = x^2+bx+c=0$ as a geometric problem 
about the area of some region which when suitably completed becomes a square. 

\section*{Appendix C: Mathematica code}

Here are some of the basic functions: 

\lstset{language=Mathematica} \lstset{frameround=fttt}
\begin{lstlisting}[frame=single]
X[x_,n_]:=Product[x-j,{j,0,n-1}]; 
sin[a_,x_]:=Im[exp[a I,x]];    cos[a_,x_]:=Re[exp[a I,x]];
tan[a_,x_]:=sin[a,x]/cos[a,x]; cot[a_,x_]:=cos[a,x]/sin[a,x]; 
exp[a_,x_]:=(1+a)^x;           log[a_,x_]:=Log[1+a,x];
rec[x_]:=log[1,x+1]-log[1,x]; 
\end{lstlisting}

Here is an illustration of the fundamental theorem of calculus: 

\lstset{language=Mathematica} \lstset{frameround=fttt}
\begin{lstlisting}[frame=single]
d[f_] := Function[m, f[m + 1] - f[m]];
s[f_] := Function[m, Sum[f[k],{k,0,m-1}]];
f0=Function[x,sin[1,x]]; g0=d[f0]; h0=s[g0]; 
{f0[201],h0[201]}
\end{lstlisting}

And here is an example of a Taylor expansion leading to en
interpolation. 

\lstset{language=Mathematica} \lstset{frameround=fttt}
\begin{lstlisting}[frame=single]
f0[x_]:=Cos[x^2/2]; M=10;
d[f_]:=Function[m,f[m+1]-f[m]];
Di[f_,k_]:=Module[{g,h=f},Do[g=d[h];h=g,{j,k}];h];
X[x_, n_]:=Product[(x-j),{j,0,n-1}]; 
Taylor[f_]:=Table[g=Di[f,k];g[0],{k,0,M}];
a=Taylor[f0]; f1[x_]:=Sum[X[x,j] a[[j+1]]/j!,{j,0,M}];
S1=Plot[{f0[x],f1[x]},{x,0,M}];
S2=Graphics[Table[{Red,Point[{k,f1[k]}]},{k,0,M-1}]];
Show[{S1,S2}]
\end{lstlisting}

The functions can also be defined using more general ``Planck 
constant" $h$: 

\lstset{language=Mathematica} \lstset{frameround=fttt}
\begin{lstlisting}[frame=single]
exp[a_,h_,x_]:= (1 + a h)^(x/h); 
cos[a_,h_,x_]:=Re[exp[I a,h,x]]; 
sin[a_,h_,x_]:=Im[exp[I a,h,x]]; 
Plot[{Cos[x],cos[1,0.1,x],Sin[x],sin[1,0.1,x]},{x,0,4Pi}]
\end{lstlisting}

Here are some graphs, already built into the computer algebra system:

\lstset{language=Mathematica} \lstset{frameround=fttt}
\begin{lstlisting}[frame=single]
CompleteGraph[5]
CycleGraph[7]
WheelGraph[8]
StarGraph[11]
PD = PolyhedronData; p=PD["Platonic"];
UG = UndirectedGraph
UG[Graph[PD[p[[3]],"SkeletonRules"]]]
UG[Graph[PD[p[[4]],"SkeletonRules"]]]
\end{lstlisting}

\begin{figure}[H]
\scalebox{0.52}{\includegraphics{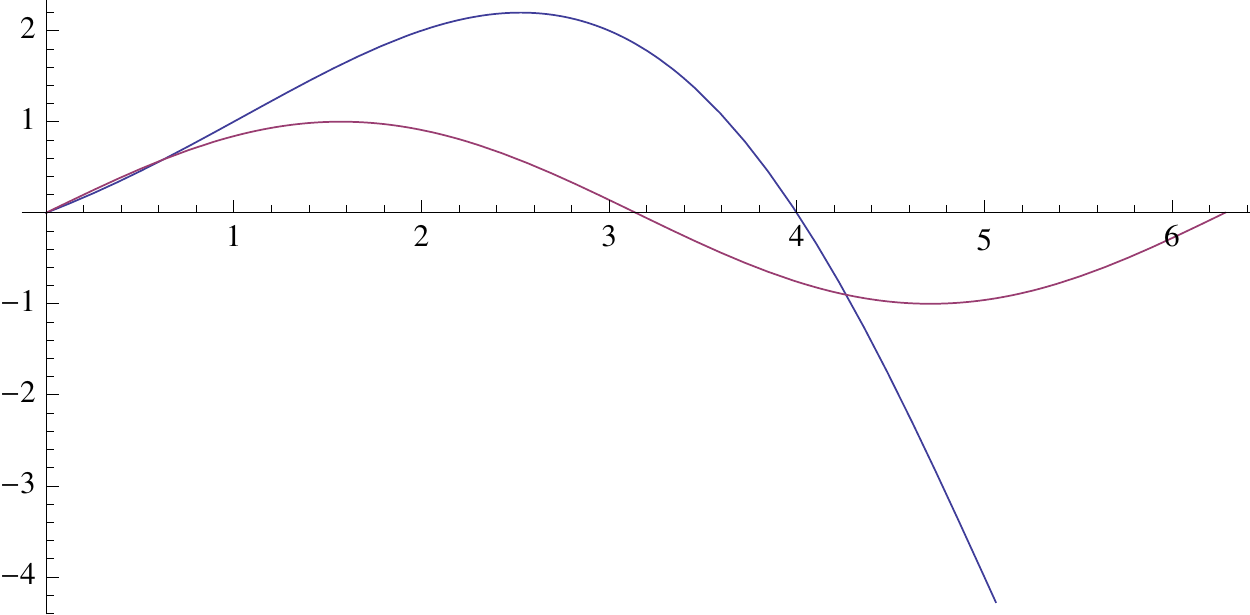}}
\scalebox{0.52}{\includegraphics{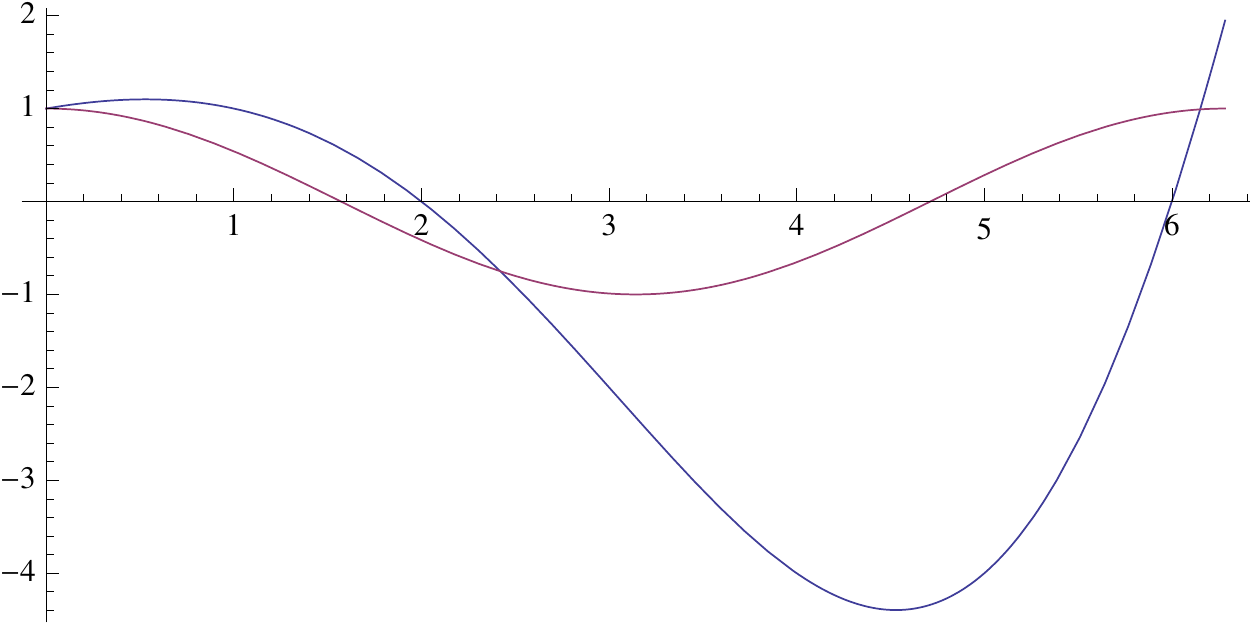}}
\scalebox{0.52}{\includegraphics{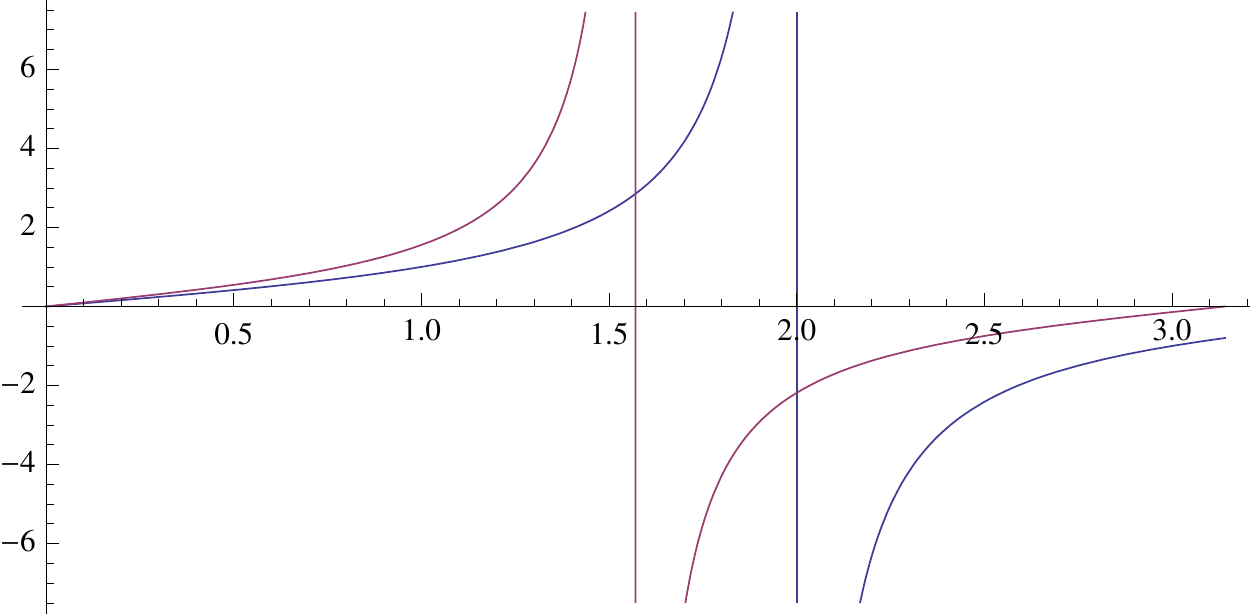}}
\scalebox{0.52}{\includegraphics{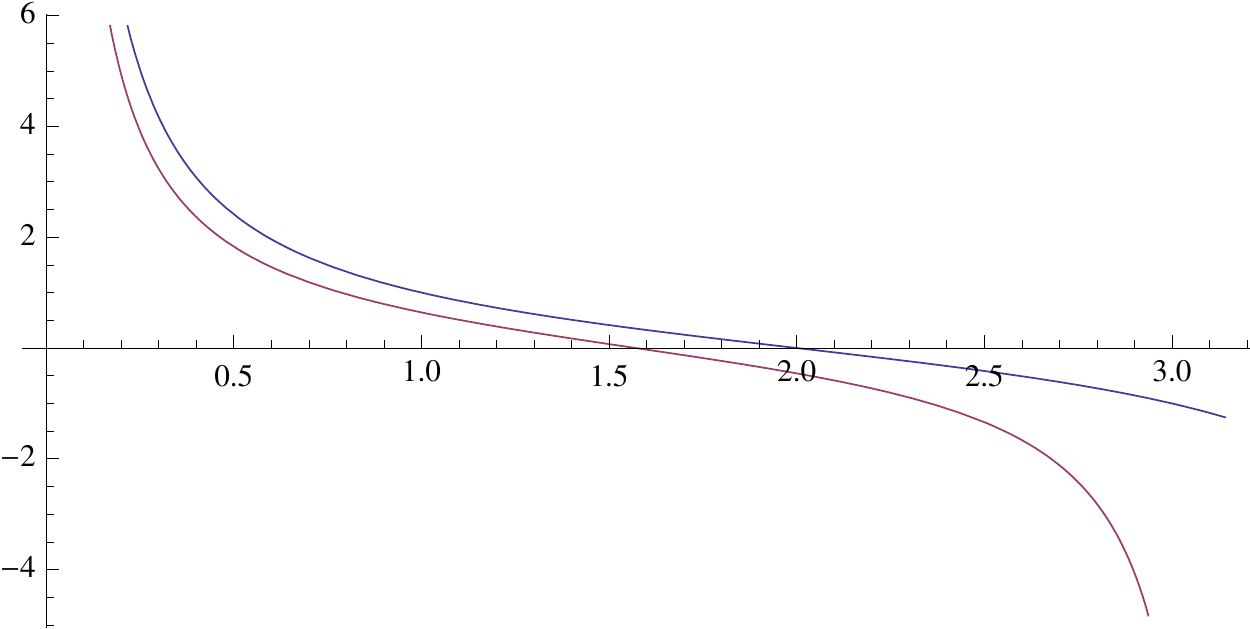}}
\scalebox{0.52}{\includegraphics{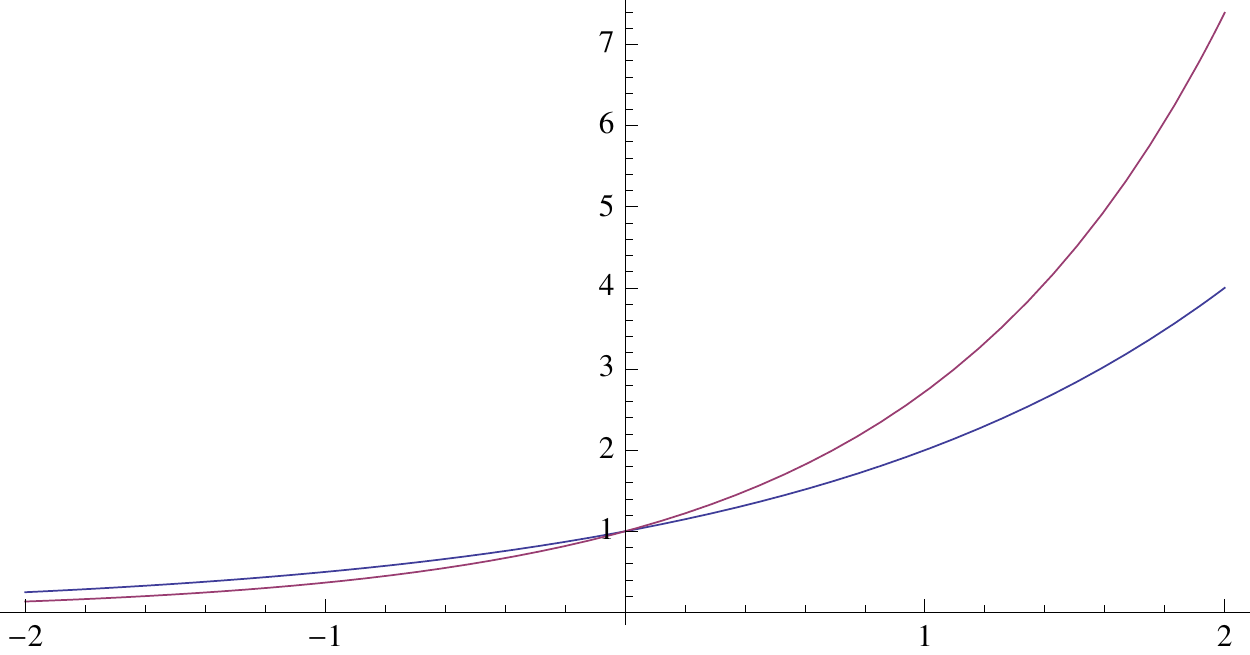}}
\scalebox{0.52}{\includegraphics{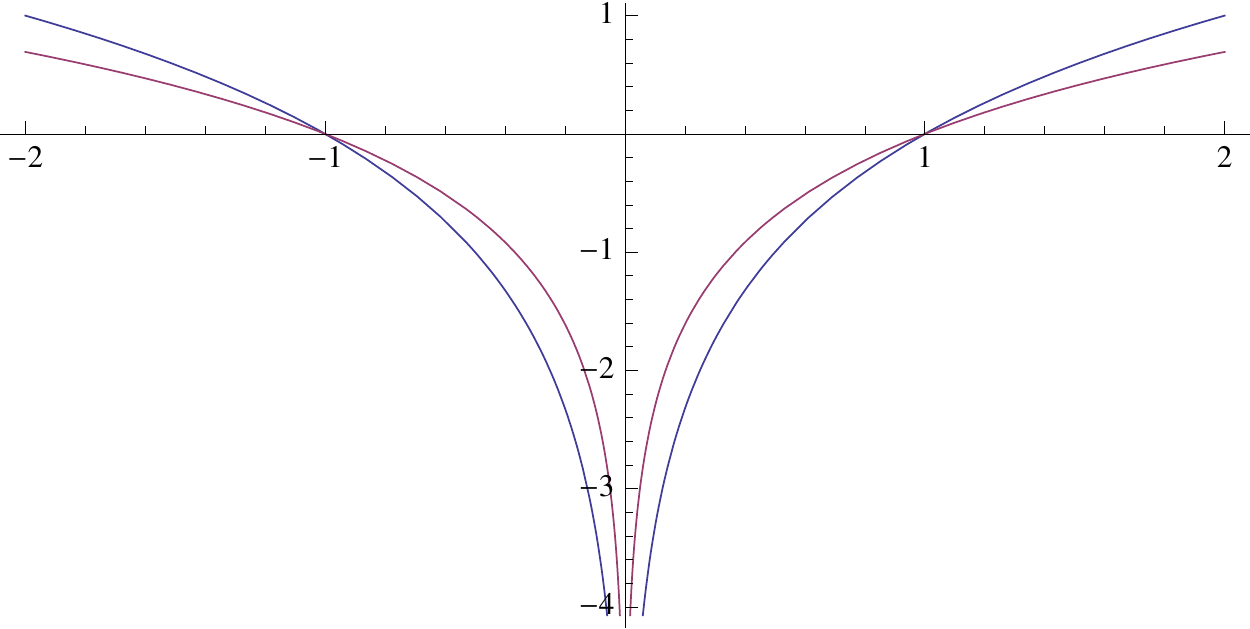}}
\scalebox{0.52}{\includegraphics{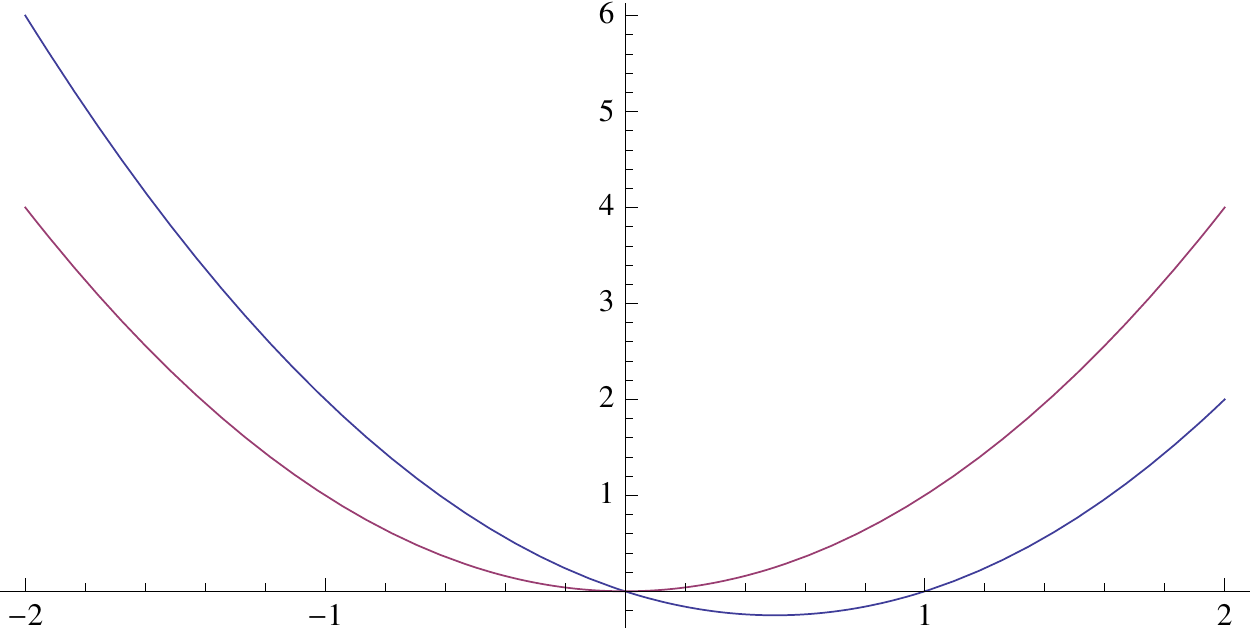}}
\scalebox{0.52}{\includegraphics{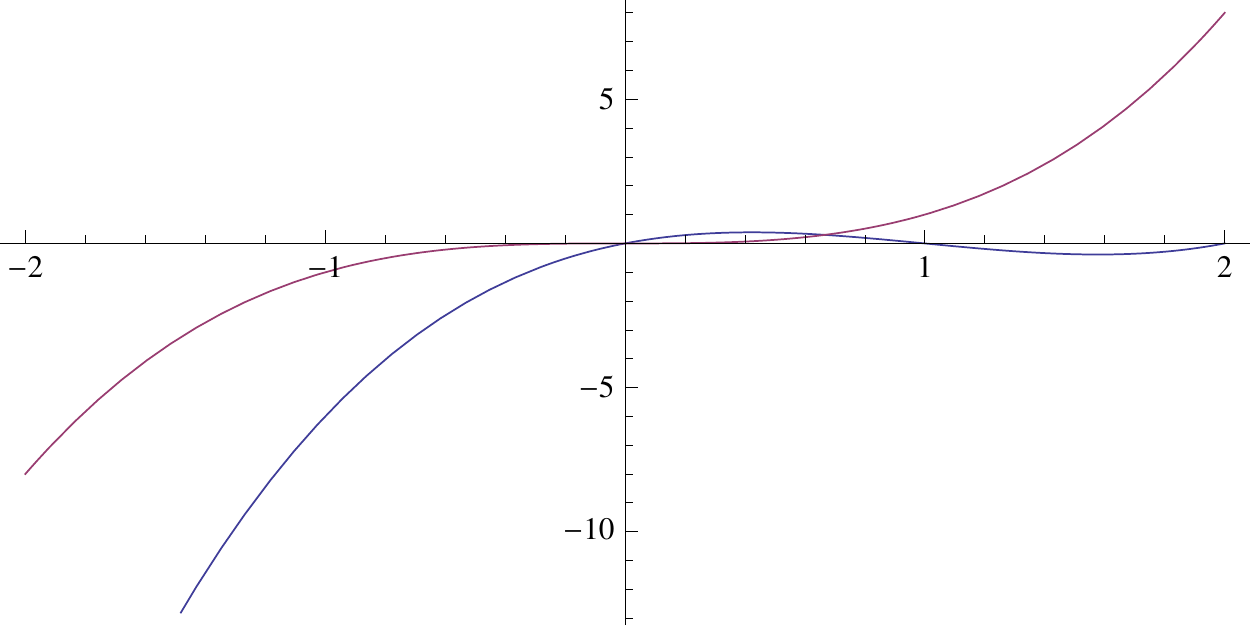}}
\caption{
Basic functions $\sin,\cos,\tan,\cot,\exp,\log,x^2,x^3$
compared with the classical functions. 
}
\end{figure}

\begin{figure}[H]
\scalebox{0.82}{\includegraphics{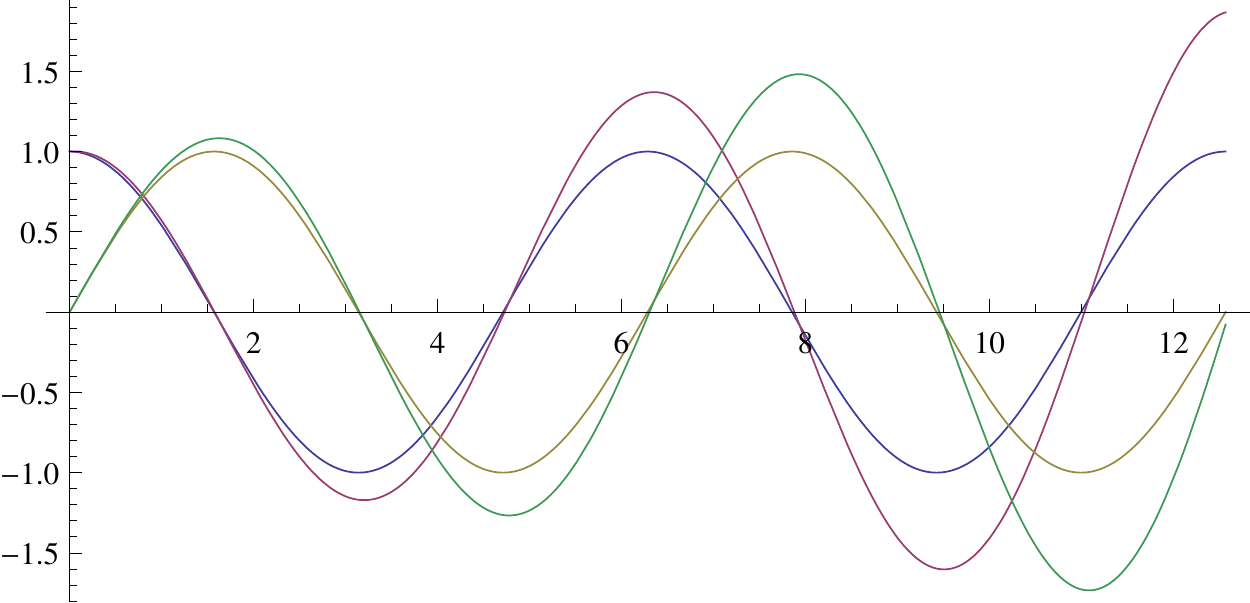}}
\caption{
A comparison of $\sin_h(x),\sin(x)$ and $\cos_h(x),\cos(x)$
for $h =0.1$. }
\end{figure}

\begin{figure}[H]
\scalebox{0.99}{\includegraphics{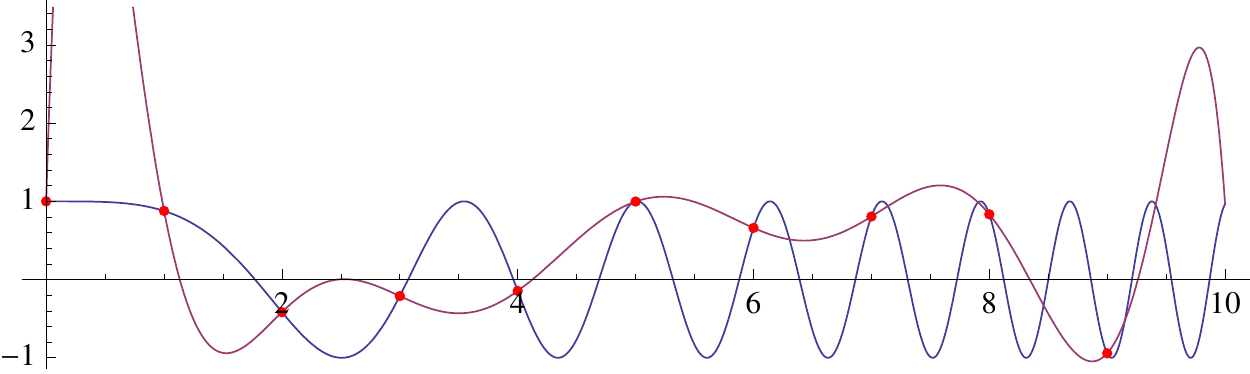}}
\caption{
Interpolating the function $\cos(x^2/2)$ using the discrete Taylor
formula of Newton and Gregory. 
}
\end{figure}

\begin{figure}[H]
\scalebox{0.38}{\includegraphics{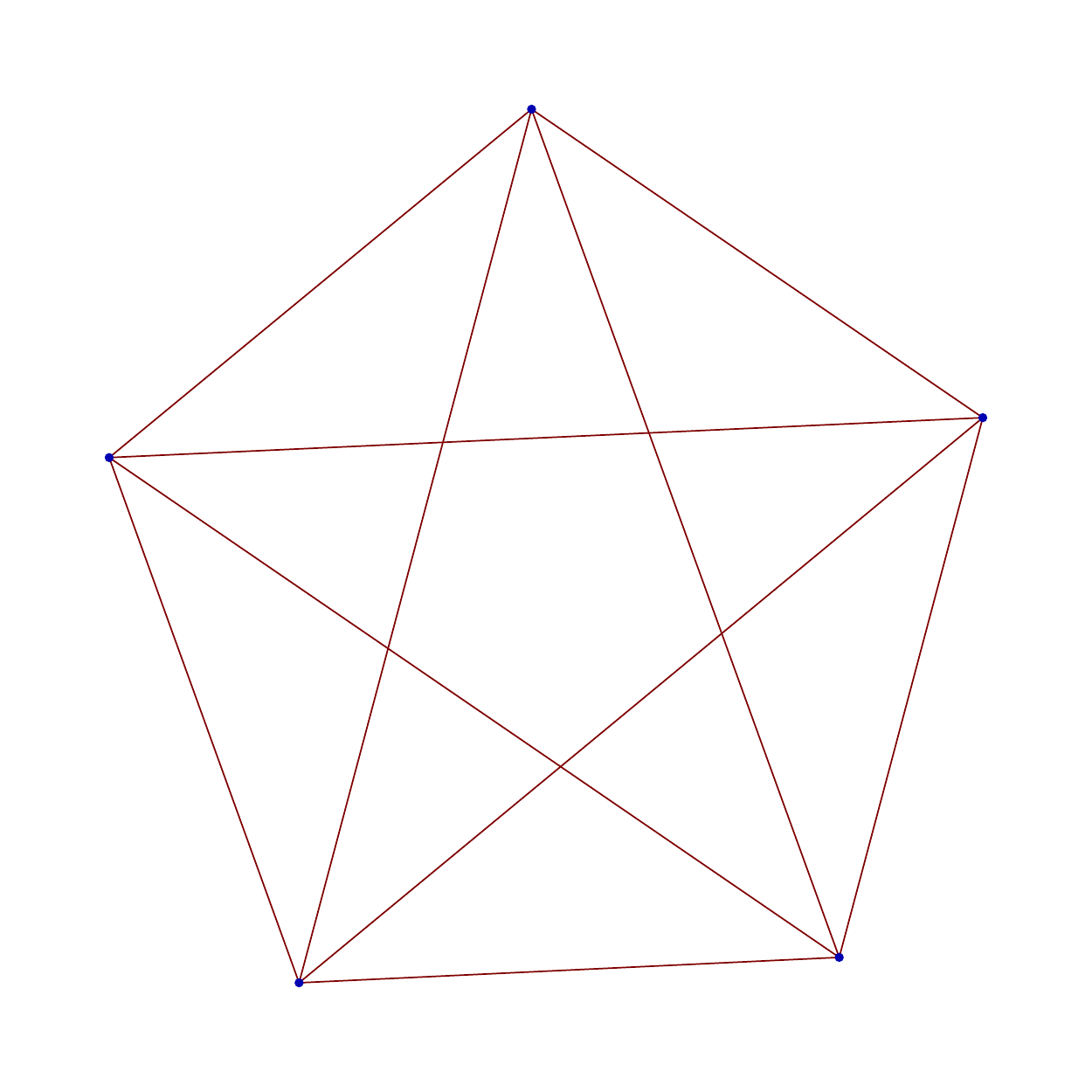}}
\scalebox{0.38}{\includegraphics{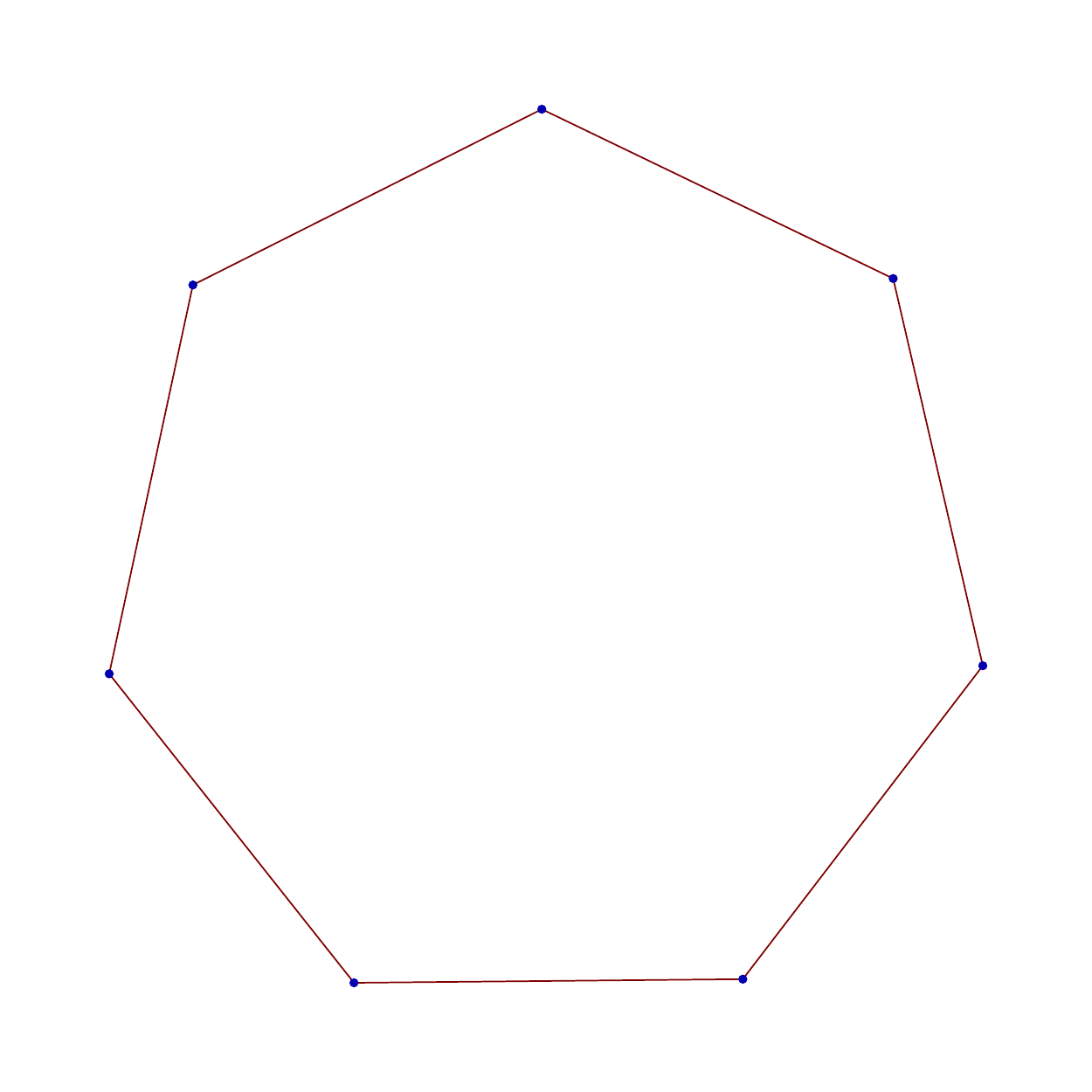}}
\scalebox{0.38}{\includegraphics{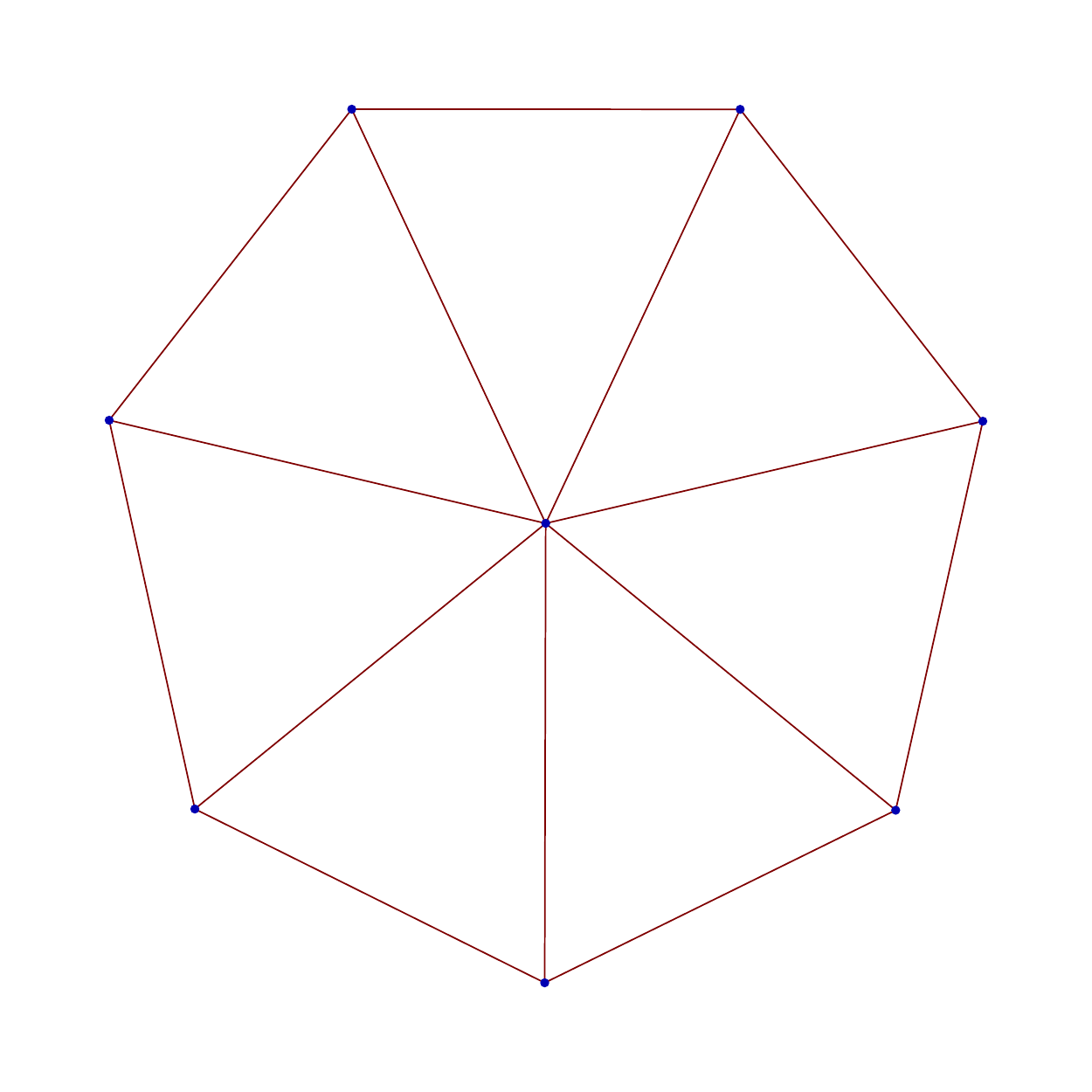}}
\scalebox{0.38}{\includegraphics{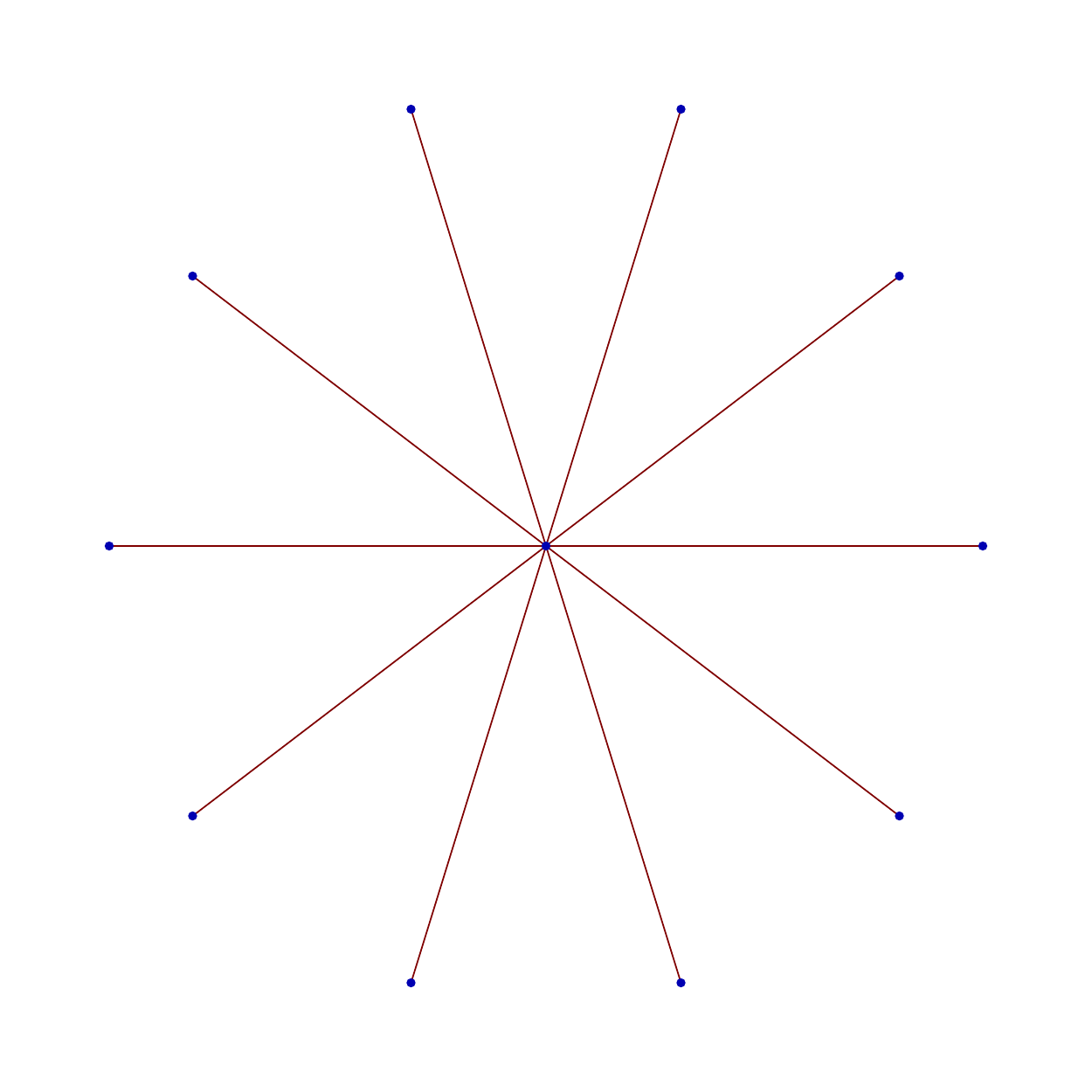}}
\scalebox{0.38}{\includegraphics{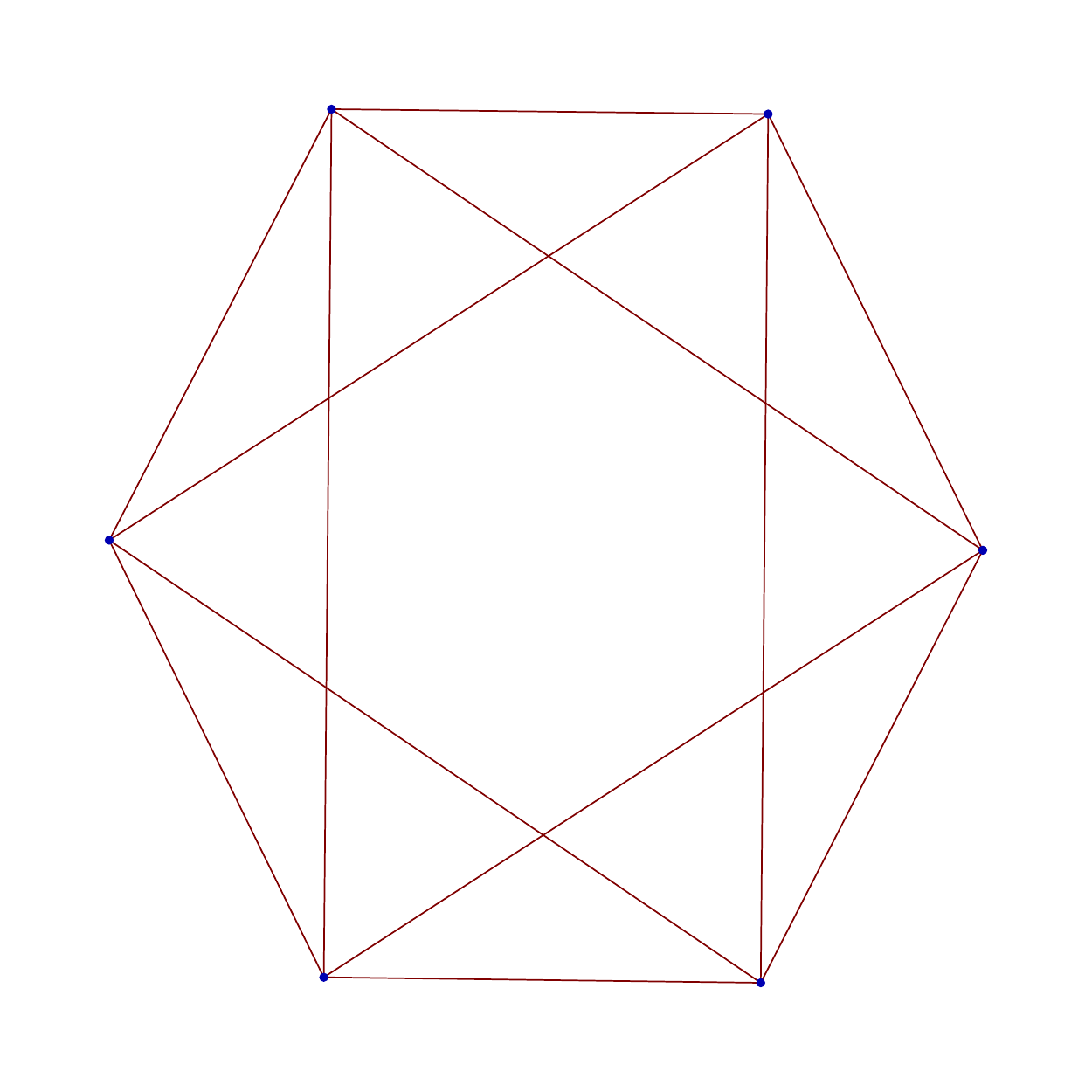}}
\scalebox{0.38}{\includegraphics{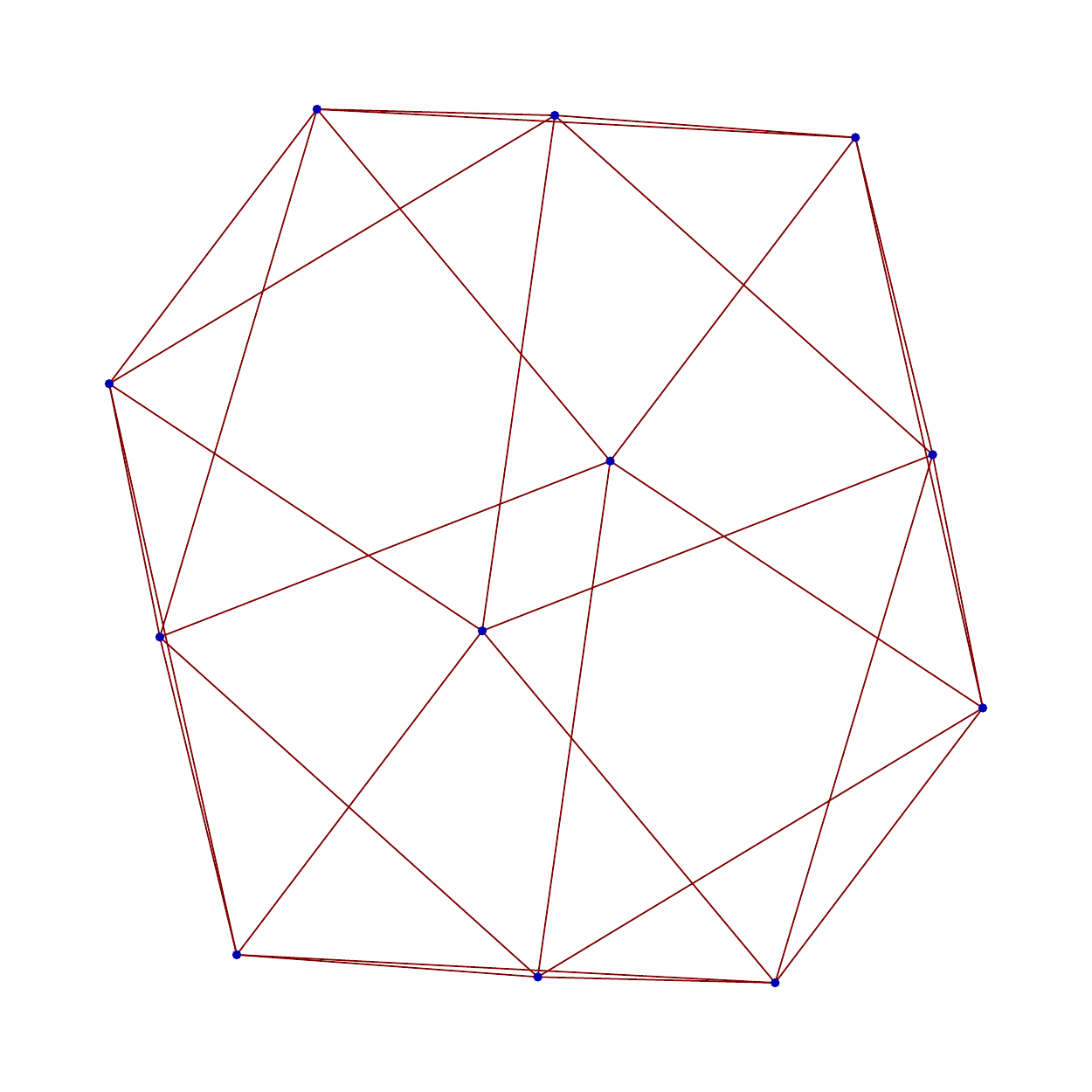}}
\caption{
Example graphs $K_5,C_7,W_7,S_{10},O,I$.
}
\end{figure}

\section{Calculus flavors}

Calculus exists in many different flavors. We have infinite-dimensional 
versions of calculus like {\bf functional analysis} and {\bf calculus of variations},
calculus has been extended to less regular geometries and functions using
the language of {\bf geometric measure theory}, {\bf integral geometry} or the {\bf theory of 
distributions}. Calculus in the complex in the form of {\bf complex analysis} has
become its own field and calculus on more general spaces is known as
{\bf differential topology}, {\bf Riemannian geometry} or {\bf algebraic geometry}, depending
on the nature of the objects and the choice of additional structures.  \\

The term ``quantum calculus" is used in many different contexts as mathematical objects can be
``quantized" in different ways. One can discretize, replace commutative algebras
with non-commutative algebras, or replace classical variational problems with path
integrals. ``Quantized calculus" is part of a more {\bf general quantum calculus},
where a measure $\mu$ on the real line is chosen leading to the derivative 
$f' =\int (f(x+h)-f(x))/h \; d\mu(h)$. 
{\bf Quantized calculus} is a special case, when $\mu$ is the Lebesgue measure and
{\bf quantum calculus} is the case when $\mu$ is a Dirac measure so that $f' = (f(x+h)-f(x))/h$. 
The measure $\mu$ defining the derivative has a smoothing effect if it is applied to space of derivatives.
\cite{Connes} defines quantized calculus as a derivative $df=[F,f]$, where $F:H\to H$ is a selfadjoint
operator on a Hilbert space $H$ satisfying $F^2=1$. The Hilbert transform has this property.
Real variables are replaced by self-adjoint operator, {\bf infinitesimals}
are {\bf compact operators} $T$, the {\bf integral} is the {\bf Dixmier trace} $||T||$ and integrable functions are replaced
by compact operators for which $||T|| = \lim (1/\log(n)) \sum_{k \leq n} \lambda_k$ exists in a suitable sense,
where $\lambda_k$ are the eigenvalues of the operator $|T|$ arranged in decreasing order. 
The Dixmier trace neglects infinitesimals, compact operators $T$ for which $\lambda_k = o(k^{-1})$. 
A {\bf Fredholm module} is a representation of an algebra $A$ as operators in the Hilbert space $H$
such that $[F,f]$ is infinitesimal (=compact) for all $f \in A$. In the simplest case, for
single variable calculus, $F$ is the Hilbert transform. \\

It is evident that quantum calculus approaches like \cite{Connes} (geared towards geometry) and 
\cite{KacCheung} (geared towards number theory) would need considerable 
adaptations to be used in a standard calculus sequence for non-mathematicians. Quantum calculus flavors
already exist in the form of {\bf business calculus}, where the concept of limit is discarded. 
For a business person, calculus can be done on spreadsheets, integration is the process of adding 
up rows or columns and differentiation is looking at trends, the goal of course being to integrate
up the changes and get trends for the future. Books like \cite{HoffmannBradley} (a ``brief edition" with 800 pages)
however follow the standard calculus methodology. \\

I personally believe that the standard calculus we teach today is and will remain an important 
benchmark theory which has a lot of advantages over discretized calculus flavors. From experience, I noticed 
for example that having to carry around an additional discretization
parameter $h$ can be a turn-off for students. This is why in the talk, $h=1$
was assumed. In the same way as {\bf classical mechanics} is an elegant idealization of 
{\bf quantum mechanics} and {\bf Newtonian mechanics} is a convenient limit of {\bf general relativity}, 
classical calculus is a teachable idealisation of a more general quantum calculus.  \\

Still, a discrete approach to calculus could be an option for
a course introducing proofs or research ideas. Many of the current
calculus books have grown to huge volumes and resemble each other even to the point
that different textbooks label chapters in the same order and that problems resemble each other,
many of them appear generated by computer algebra systems.  \\

But as many newcomers have seen, the transition from a standard textbook to a new textbook is hard. 
If there are too many changes even if it only affects the order in which the material is taught, 
then this can alienate the customer. Too fast innovation destroys its sales. 
Innovative approaches are seldom successful commercially, especially because college
calculus systems have become complicated machineries, where larger courses with 
different sections are coordinated and teams of teachers 
and course assistants work together. When looking at the history of calculus
textbooks from the past, over decades or centuries, we can see that the development
today is still astonishingly fast today and that the variety has never been bigger.   \\

A smaller inquiry based or proof based course aiming to have students
in a ``state of research" can benefit from fresh approaches. The reason
is that asking to reprove standard material have become 
almost impossible today. It is due to the fact that one can look up proofs so 
easily today: just enter the theorem into a search engine and almost certainly a 
proof can be found. It is no accident 
that more inquiry based learning stiles have flourished in a time when point set topology 
has been developed. Early point set topology for example was a perfect toy to 
sharpen research skills. It has become almost impossible now to 
make any progress in basic set topology because the topic has been cataloged 
so well. When developing a new theorem, a contemporary student almost certainly 
will recover a result which is a special case of something which has been 
found 50-100 years ago. In calculus it is even more extreme because the theorems
taught have been found 200-300 years ago.  \\

Learning old material in a new form is always exciting. I myself have learned
calculus the standard way but then seen it again in a ``nonstandard analysis course" 
given by Hans L\"auchli. The course had been based on the article \cite{Nelson77}. 
We were also recommended to read \cite{Robert}, a book which stresses 
in the introduction how Euler has worked already in a discrete setting. Euler would
write $e^x = (1+x/n)^n$, where $n$ is a number of infinite size 
(Euler: "n numerus infinite magnus") (\cite{Euler} paragraph 155). 
Nelson also showed that probability theory on finite sets works well \cite{Nelson}. 
There are various textbooks which tried nonstandard analysis. One which is influenced
by the Robinson approach to nonstandard analysis is \cite{Keisler}. 
But nonstandard analysis has never become the standard. It remains true to its name. 

\section{Calculus textbooks}

The history of calculus textbooks is fascinating and would deserve to be written down
in detail. Here are some comments about some textbooks from my own calculus library.
The first textbook in calculus is L'H\^opital's 
``Analyse des Infiniment Petits pour l'Intelligence des Lignes Courbes"  \cite{lHopital},
a book from 1696, which had less than 200 pages and still contains 156 figures. 
Euler's textbook \cite{Euler} on calculus starts the text with a modern treatment of
functions but still emphasis "analytic expressions". Euler for example would hardly 
considered the rule which assigns $1$ to an irrational numbers and $0$ on rational number,
to be a function. Euler uses already in that textbook already his formula 
$e^{ix} = \cos(x) + i \sin(x)$ extensively but in material is closer of what we
consider ``precalculus" today, even so it goes pretty far, especially with
series and continued fraction expansions. An other master piece is
the textbook of Lagrange \cite{Lagrange}. 
\cite{KatzHistoria} mentions textbooks of the 18th century with authors 
C. Reyneau,G. Chayne,C. Hayes,H. Ditton,J. Craig,E. Stone,J. Hodgson, J. Muller or T. Simpson. 
Hardy's book \cite{Hardy}, whose first edition came out in 1908 also contains a 
substantial real analysis part and uses already the modern function approach. It no more
restricts to ``analytic expressions" but also allows functions to be limiting
expressions of series of functions. It is otherwise a single variable text book. 
From about the same time is \cite{Murray}, which is more than 100 years old but resembles
in style and size modern textbooks already, including lots of problems, also some of 
applied nature. 
The book \cite{GibbsWilson} is a multivariable text which has many applications to 
physics and also contains some differential geometry. 
The book \cite{Laisant} is remarkable because it might be where we ``stole" the first 
slide of our talk from, mentioning sticks and pebbles to make the first steps in calculus. 
This book weights only 150 pages but covers lots of ground. 
\cite{Osgood} treats single variable on 200 pages with a decent amount of 
exercises and \cite{Hay} manages with the same size to cover a fair amount of multivariable
calculus and tensor analysis, differential geometry and calculus of variations
in three dimensions, uses illustrations however only sparingly. The structure of Hay's book 
is in parts close to modern calculus textbooks like \cite{OstebeeZorn,Stewart,Thomas,Anton}.
many of them started to grow fatter. The textbook \cite{Widder} for example has already 420
pages but covers single and multivariable calculus, 
real analysis and measure theory, Fourier series and partial differential equations. But it is
no match in size to books like \cite{Anton} with 1318 pages or Stewart with about the same size.
Important in Europe was \cite{Courant} which had been translated into other
language and which my own teachers still recommended. Its influence is clearly visible
in all major textbooks which are in use today. \cite{Lang} made a big
impact as it is written with great clarity.  \\

Since we looked at a discrete setting, we should mention that calculus books based
on differences have appeared already 60 years ago. The book \cite{Jordan} illustrates
the difficulty to find good notation. This is a common obstacle with discrete
approaches, especially also with numerical analysis books. Discrete calculus often suffers
from ugly notation using lots of indices.  One of the reason why Leibnitz calculus has 
prevailed is the use of intuitive notation. For a modern approach see \cite{GradyPolimeni}. \\

Not all books exploded in size. \cite{Gleason} covers in 120 pages 
quite a bit of real analysis and multivariable calculus. 
But books would become heavier. \cite{Agnew} already counted 725 pages. It was intended for
a three semester course and stated in the introduction: {\it there is an element of 
truth in the old saying that the Euler textbook 'Introductio in Analysin Infinitorum,
(Lausannae, 1748)' was the first great calculus textbook and that all elementary 
calculus textbooks copied from books that were copied from Euler.} 
The classic \cite{Kline} had already close to 1000 pages and already
included a huge collection of problems. 
Influential for the textbooks of our modern time are also Apostol \cite{Apostol1,Apostol2}
which I had used myself to teach. It uses elements of linear algebra, 
probability theory, real analysis, complex numbers and differential equations. Its
influence is also clearly imprinted on other textbooks. The book \cite{Edwards} was 
one of the first which dared to introduce calculus in arbitrary dimensions using differential forms. 
Since then, calculus books have diversified and exploded in general in size. 
There are also exceptions: original approaches include 
\cite{DivGradCurl,Sparks,Nahin1,Levi,Niven,AdamsThompsonHass,Manga,Ekeland,Gonick,Paradoxes}. 

\section{Tweets}

We live in an impatient twitter time. Here are two actually tweeted 140 character
statements which together cover the core of the fundamental theorem of calculus:  \\

\begin{center}
\fbox{ \parbox{10cm}{ 
With $s(x)=x+1$ \\
and inverse $t(x)=x-1$, \\
define $X^n=(x t)^n$. \\
For example, \\ 
$1=1$, $X=x$, \\
$X^2=x(x-1)$, \\
$X^3=x(x-1)(x-2)$. 
}}
\end{center}

and  \\

\begin{center}
\fbox{ \parbox{10cm}{ 
$Df(x)=[s,f]=f(x+1)-f(x)$ \\
gives $D X^n=n X^{n-1}$. With \\
$Sf(x)=f(0)+f(1)+...+f(n-1)$, \\
the fundamental theorem \\
$SDf(x)=f(x)-f(0)$  \\
$DS f(x)=f(x)$ \\
holds.

}}
\end{center}


\begin{thebibliography}{10}

\bibitem{AdamsThompsonHass}
C.~Adams, A.~Thompson, and J.~Hass.
\newblock {\em How to Ace Calculus, The Streetwise guide}.
\newblock Henry Holt and Company, LLC, 1998.

\bibitem{Agnew}
R.P. Agnew.
\newblock {\em Calculus, Analytic Geometry and Calculus, with Vectors}.
\newblock McGraw-Hill Book Company, Inc, New York, 1962.

\bibitem{Apostol1}
T.M. Apostol.
\newblock {\em One-variable calculus, introduction to linear algebra}.
\newblock John Wiley and Suns, New York, second edition edition, 1967.

\bibitem{Apostol2}
T.M. Apostol.
\newblock {\em Calculus: Multivariable Calculus and Linear algebra with
  applications to differential equations and probability theory}.
\newblock John Wiley and Suns, New York, second edition edition, 1969.

\bibitem{Arnold}
V.I. Arnold.
\newblock {\em Huygens and Barrow, Newton and Hooke}.
\newblock Birkh{\"a}user, 1990.

\bibitem{Bardi}
J.S. Bardi.
\newblock {\em The Calculus wars}.
\newblock Thunder's Mouth Press, New York, 2006.

\bibitem{Bell}
E.T. Bell.
\newblock {\em Men of Mathematics, the lives and achievements of the great
  mathematicians from Zeno to Poincare}.
\newblock Penguin Books, Melbourne, London, Baltimore, 1937-1953.

\bibitem{Bottazzini}
U.~Bottazzini.
\newblock {\em The Higher Calculus, A History of Real and Complex Analysis from
  Euler to Weierstrass}.
\newblock Springer Verlag, 1986.

\bibitem{BoyerCalculus}
C.B. Boyer.
\newblock {\em The History of the Calculus and its Conceptual Development}.
\newblock Dover Publications, 2 edition, 1949.

\bibitem{Bressoud}
D.~M. Bressoud.
\newblock Historical reflections on teaching the fundamental theorem of
  integral calculus.
\newblock {\em American Mathematical Monthly}, 118:99--115, 2011.

\bibitem{Cajori}
F.~Cajori.
\newblock {\em A history of Mathematics}.
\newblock MacMillan and Co, 1894.

\bibitem{Connes}
A.~Connes.
\newblock {\em Noncommutative geometry}.
\newblock Academic Press, 1994.

\bibitem{Courant}
R.~Courant.
\newblock {\em Differentrial and Integral Calculus, Volume I and II}.
\newblock Interscience Publishers, Inc, New York, 1934.

\bibitem{lHopital}
G.~de~L'{H\^opital}.
\newblock {\em Analyse des Infiniment Petits pour l'Intelligence des Lignes
  Courbes}.
\newblock Francois Montalant, Paris, second edition, 1716.
\newblock first edition 1696.

\bibitem{DieudonneMusic}
J.~A. Dieudonne.
\newblock {\em Mathematics - the music of reason}.
\newblock Springer Verlag, 1998.

\bibitem{Dunham}
W.~Dunham.
\newblock {\em The Calculus Gallery, Masterpieces from Newton To Lebesgue}.
\newblock Princeton University Press, 2005.

\bibitem{DunhamMonthly}
W.~Dunham.
\newblock Touring the calculus gallery.
\newblock {\em Amer. Math. Monthly}, 112(1):1--19, 2005.

\bibitem{EdwardsHistory}
C.H. Edwards.
\newblock {\em The historical Development of the Calculus}.
\newblock Springer Verlag, 1979.

\bibitem{Edwards}
H.M. Edwards.
\newblock {\em Advanced Calculus, A differential Forms Approach}.
\newblock Birkh{\"a}user, 1969-1994.

\bibitem{Ekeland}
I.~Ekeland.
\newblock {\em Exterior differnential Calculus and Applications to Economic
  Theory}.
\newblock Sculoa Normale Superiore, Pisa, 1998.

\bibitem{Euler}
L.~Euler.
\newblock {\em Einleitung in die Analysis des Unendlichen}.
\newblock Julius Springer, Berlin, 1885.

\bibitem{FourierChaleur}
J-B.J. Fourier.
\newblock {\em Th{\'e}orie Analytique de la Chaleur}.
\newblock Cambridge University Press, 1822.

\bibitem{GibbsWilson}
J.W. Gibbs and E.B. Wilson.
\newblock {\em Vector Analysis}.
\newblock Yale University Press, New Haven, 1901.

\bibitem{Gleason}
A.~M. Gleason.
\newblock {\em Linear analysis and Calculus, Mathematics 21}.
\newblock Harvard University Lecture Notes, 1972-1973.

\bibitem{Gonick}
L.~Gonick.
\newblock {\em The Cartoon Guide to Calculus}.
\newblock William Morrow, 2012.

\bibitem{GradyPolimeni}
L.J. Grady and J.R. Polimeni.
\newblock {\em Discrete Calculus, Applied Analysis on Graphs for Computational
  Science}.
\newblock Springer Verlag, 2010.

\bibitem{Anton}
I.C.~Bivens H.~Anton and S.~Davis.
\newblock {\em Calculus, 10th edition}.
\newblock Wiley, 2012.

\bibitem{Hardy}
G.H. Hardy.
\newblock {\em A course of pure mathematics}.
\newblock Cambridge at the University Press, third edition, 1921.

\bibitem{Hay}
G.E. Hay.
\newblock {\em Vector and Tensor analysis}.
\newblock Dover Publications, 1953.

\bibitem{HoffmannBradley}
L.D. Hoffmann and G.L. Bradley.
\newblock {\em Calculus, for Business, Economics, and the social life
  sciences}.
\newblock McGraw Hill, brief edition, 2007.

\bibitem{Jordan}
C.~Jordan.
\newblock {\em Calculus of finite differences}.
\newblock Chelsea publishing company, New York, 1950.

\bibitem{KacCheung}
V.~Kac and P.~Cheung.
\newblock {\em Quantum calculus}.
\newblock Universitext. Springer-Verlag, New York, 2002.

\bibitem{KatzHistoria}
V.J. Katz.
\newblock The calculus of the trigonometric functions.
\newblock {\em Historia Math.}, 14(4):311--324, 1987.

\bibitem{Keisler}
H.J. Keisler.
\newblock {\em Elementary Calculus, An infinitesimal approach}.
\newblock Creative Commons, second edition, 2005.

\bibitem{Kline}
M.~Kline.
\newblock {\em Calculus}.
\newblock John Wiley and Sons, 1967-1998.

\bibitem{Paradoxes}
S.~Klymchuk and S.~Staples.
\newblock {\em Paradoxes and Sophisms in Calculus}.
\newblock MAA, 2013.

\bibitem{cherngaussbonnet}
O.~Knill.
\newblock A graph theoretical {Gauss-Bonnet-Chern} theorem.
\newblock {\\}http://arxiv.org/abs/1111.5395, 2011.

\bibitem{poincarehopf}
O.~Knill.
\newblock A graph theoretical {Poincar\'e-Hopf} theorem.
\newblock {\\} http://arxiv.org/abs/1201.1162, 2012.

\bibitem{indexexpectation}
O.~Knill.
\newblock On index expectation and curvature for networks.
\newblock {\\}http://arxiv.org/abs/1202.4514, 2012.

\bibitem{knillcalculus}
O.~Knill.
\newblock {The theorems of Green-Stokes,Gauss-Bonnet and Poincare-Hopf in Graph
  Theory}.
\newblock {\\}http://arxiv.org/abs/1201.6049, 2012.

\bibitem{DiracKnill}
O.~Knill.
\newblock The {D}irac operator of a graph.
\newblock {\\}http://http://arxiv.org/abs/1306.2166, 2013.

\bibitem{classicalstructures}
O.~Knill.
\newblock Classical mathematical structures within topological graph theory.
\newblock {\\}http://arxiv.org/abs/1402.2029, 2014.

\bibitem{KnillHerran}
O.~Knill and O.~Ramirez-Herran.
\newblock Space and camera path reconstruction for omni-directional vision.
\newblock Retrieved February 2, 2009, from http://arxiv.org/abs/0708.2438,
  2007.

\bibitem{Manga}
H.~Kojama.
\newblock {\em Manga Guide to Calculus}.
\newblock Ohmsha, 2009.

\bibitem{Lagrange}
J.L. Lagrange.
\newblock {\em Lecons sur le calcul des fonctions}.
\newblock Courcier, Paris, 1806.

\bibitem{Laisant}
C.A. Laisant.
\newblock {\em Mathematics, Illustrated}.
\newblock Constable and Company Ltd, London, 1913.

\bibitem{Lang}
S.~Lang.
\newblock {\em A first course in Calculus}.
\newblock Springer, fifth edition, 1978-1986.

\bibitem{Levi}
M.~Levi.
\newblock {\em The Mathematical Mechanic}.
\newblock Princeton University Press, 2009.

\bibitem{Luzin}
N.~Luzin.
\newblock Function. {I},{II}.
\newblock {\em Amer. Math. Monthly}, 105(1):59--67, 263--270, 1998.

\bibitem{Murray}
D.A. Murray.
\newblock {\em A first course in infinitesimal calculus}.
\newblock Longmans, Green and Co, New York, 1903.

\bibitem{Nahin1}
Paul~J. Nain.
\newblock {\em When least is best}.
\newblock Princeton University Press, 2004.

\bibitem{Nelson77}
E.~Nelson.
\newblock Internal set theory: A new approach to nonstandard analysis.
\newblock {\em Bull. Amer. Math. Soc}, 83:1165--1198, 1977.

\bibitem{Nelson}
E.~Nelson.
\newblock {\em Radically elementary probability theory}.
\newblock Princeton university text, 1987.

\bibitem{Niven}
I.~Niven.
\newblock {\em Maxima and Minima without Calculus}, volume~6 of {\em Dolciani
  Mathematical Expositions}.
\newblock MAA, 1981.

\bibitem{Osgood}
W.F. Osgood.
\newblock {\em Elementary Calculus}.
\newblock MacMillan Company, 1921.

\bibitem{OstebeeZorn}
A.~Ostebee and P.~Zorn.
\newblock {\em Multivariable Claculus}.
\newblock Sounders College Publishing, 1998.

\bibitem{Robert}
A.~Robert.
\newblock {\em Analyse non standard}.
\newblock Presses polytechniques romandes, 1985.

\bibitem{Rosenthal}
A.~Rosenthal.
\newblock The history of calculus.
\newblock {\em Amer. Math. Monthly}, 58:75--86, 1951.

\bibitem{DivGradCurl}
H.M. Schey.
\newblock {\em div grad curl and all that}.
\newblock Norton and Company, New York, 1973-1997.

\bibitem{Sparks}
J.C. Sparks.
\newblock {\em Calculus without limits - Almost}.
\newblock Author House, third edition, 2007.

\bibitem{Stewart}
J.~Stewart.
\newblock {\em Calculus, Concepts and contexts}.
\newblock Brooks/Cole, 4e edition, 2010.

\bibitem{Tent}
M.B.W. Tent.
\newblock {\em Gottfried Wilhelm Leibniz, the polymath who brought us
  calculus}.
\newblock CRC Press, 2012.

\bibitem{Thomas}
G.B. Thomas, M.D. Weir, and J.R. Hass.
\newblock {\em Thomas' Calculus, Single and Multivariable, 12th edition}.
\newblock Pearson, 2009.

\bibitem{Waestlund}
J.~{W\"astlund}.
\newblock A yet simpler proof of the chain rule.
\newblock {\em American Mathematical Monthly}, 120:900, 2013.

\bibitem{Widder}
D.V. Widder.
\newblock {\em Advanced Calculus}.
\newblock Prentice-Hall, Inc, New York, 1947.

\end{thebibliography}

\end{document}